\documentclass[aos,MSNbibl,nameyear,dvips]{arximspdf}
\usepackage{graphicx}

\doi{10.1214/12-AOS1071} 
\volume{41}
\issue{1}
\pubyear{2013}
\firstpage{221}
\lastpage{249}

\makeatletter

\newcommand{\cal}{\mathcal}

\newtheorem{proposition}{Proposition}
\newtheorem{lemma}{Lemma}
\newtheorem{corollary}{Corollary}

\newproclaim{definition}{Definition}
\newtheorem{theorem}{Theorem}

\def\var{\operatorname{var}}
\def\cov{\operatorname{cov}}
\newcommand{\indep}{\perp\hspace*{-6.2pt}\perp}

\def\spc{{\cal S}}

\def\G{{\cal G}}
\def\implies{\Rightarrow}

\def\S{{\mathfrak{S}}}

\def\H{{\cal{H}}}
\def\coxx{\Sigma_{XX}}
\def\coyx{\Sigma_{YX}}
\def\coyy{\Sigma_{YY}}

\def\rxy{R_{XY}}
\def\ryx{R_{YX}}

\def\hy{{\H_{Y}}}

\def\lcentral{\spc_{ Y|X}}

\def\M{\mathfrak M}
\def\pxy{P_{ XY}}
\def\px{P_{ X}}
\def\py{P_{ Y}}
\def\borelx{{\cal F}_{ X}}
\def\borely{{\cal F}_{ Y}}

\def\pygx{P_{ Y|X}}
\def\fieldx{{\cal F}_{ X}}
\def\fieldy{{\cal F}_{ Y}}
\def\om{\Omega}

\def\omx{\Omega_{ X}}

\def\omxy{\Omega_{ XY}}
\def\fieldxy{{\cal F}_{ XY}}
\def\gfield{{\cal G}}

\def\field{{\cal F}}
\def\omy{\Omega_{ Y}}
\def\qmathbb{\mathbb Q}
\def\pmathbb{\mathbb P}
\def\sample{\Omega}

\def\R{{\cal R}}

\def\central{\S_{ Y|X}}
\def\cfield{{\G_{ Y|X}}}

\def\hxquo{{{\cal H}_{ X} }}
\def\hyquo{{{\cal H}_{ Y} }}

\def\mxx{{M_{XX}}}
\def\myy{{M_{YY}}}

\def\hone{{\cal H}_1}

\def\htwo{{\cal H}_2}

\def\ran{\operatorname{ran}}

\def\ecoyx{\tilde\Sigma_{ YX}}

\def\ecoxx{\tilde\Sigma_{ XX}}

\def\ox{\Omega_{ X}}
\def\oy{\Omega_{ Y}}
\def\real{\mathbb{R}}
\def\sp{\operatorname{span}}

\def\B{{\cal B}}
\def\C{{\cal C}}
\def\D{{\cal D}}
\def\ecoyy{\tilde{\Sigma}_{ YY}}
\def\exy{E_{ X|Y}}

\def\kappax{\kappa_{ X}}
\def\kx{K_{ X}}

\def\overran{\operatorname{\overline{ran}}}
\def\complete{{\mathfrak{C}_{ Y|X}}}

\def\vxy{V_{ X|Y}}
\def\cxy{\Sigma_{ X|Y}}

\def\cv{\operatorname{CV}}
\def\epsx{\epsilon_{ X}}
\def\epsy{\epsilon_{ Y}}
\def\gamx{\gamma_{ X}}
\def\gamy{\gamma_{ Y}}

\def\elly{L_{ Y}}

\def\diag{\operatorname{diag}}

\def\kax{\kappax}

\def\gx{G_{ X}}
\def\gy{G_{ Y}}
\def\gfam{{\cal G}_{ X}}
\newcommand{\sub}[2]{#1_{ #2}}
\def\bx{{\cal B}_{ X}}
\def\hx{{\H_{ X}}}
\def\hbasis{h_{ X}}

\def\pnx{P_{ n,X}}
\def\kym{K_{ Y}}
\def\kay{\kappa_{ Y}}
\def\cy{C_{ Y}}
\def\llx{L_{ X}}
\def\lly{L_{ Y}}
\def\llly{\ell_{ Y}}
\def\lllx{\ell_{ X}}

\newcommand{\gclass}{{\cal G}_{ X}}

\makeatother

\begin{document}
\begin{frontmatter}

\title{A general theory for nonlinear sufficient dimension reduction: Formulation and estimation}
\runtitle{Nonlinear sufficient dimension reduction}

\begin{aug}
\author[A]{\fnms{Kuang-Yao} \snm{Lee}\ead[label=e1]{kuang-yao.lee@yale.edu}},
\author[B]{\fnms{Bing} \snm{Li}\corref{}\thanksref{t1}\ead[label=e2]{bing@stat.psu.edu}}
\and
\author[B]{\fnms{Francesca} \snm{Chiaromonte}\ead[label=e3]{chiaro@stat.psu.edu}}
\runauthor{K.-Y. Lee, B. Li and F. Chiaromonte}
\affiliation{Yale University, Pennsylvania State University and Pennsylvania~State~University}
\address[A]{K.-Y. Lee\\
Department of Biostatistics\\
Yale School of Public Health\\
60 College Street\\
P.O. Box 208034\\
New Haven, Connecticut\\
USA\\
\printead{e1}} 
\address[B]{B. Li\\
F. Chiaromonte\\
Department of Statistics\\
Pennsylvania State University\\
326 Thomas Building\\
University Park, Pennsylvania 16802\\
USA\\
\printead{e2}\\
\hphantom{E-mail: }\printead*{e3}}
\end{aug}

\thankstext{t1}{Supported in part by NSF Grants DMS-08-06058 and DMS-11-06815.}

\received{\smonth{6} \syear{2012}}

%
\begin{abstract}
In this paper we introduce a general theory for nonlinear sufficient
dimension reduction, and explore its ramifications and scope. This
theory subsumes recent work employing reproducing kernel Hilbert
spaces, and reveals many parallels between linear and nonlinear
sufficient dimension reduction. Using these parallels we analyze the
properties of existing methods and develop new ones. We begin by
characterizing dimension reduction at the general level of
$\sigma$-fields and proceed to that of classes of functions, leading
to the notions of sufficient, complete and central dimension reduction
classes. We show that, when it exists, the complete and sufficient
class coincides with the central class, and can be unbiasedly and
exhaustively estimated by a generalized sliced inverse regression
estimator (GSIR). When completeness does not hold, this estimator
captures only part of the central class. However, in these cases we
show that a generalized sliced average variance estimator (GSAVE) can
capture a larger portion of the class. Both estimators require no
numerical optimization because they can be computed by spectral
decomposition of linear operators. Finally, we compare our estimators
with existing methods by simulation and on actual data sets.
\end{abstract}

%

\end{frontmatter}

\section{Introduction}\label{sec1}

In this paper we propose a general theory for nonlinear sufficient dimension
reduction (SDR), develop novel estimators and investigate their
properties under this theory.
Along with these developments we also introduce a new conditional
variance operator, which can potentially be used to
generalize all second-order dimension reduction methods to the
nonlinear case.

In its classical form, linear SDR seeks a low-dimensional linear
predictor that captures in full a
regression relationship. Imagining a regression setting that comprises
multiple predictor variables
and multiple responses, let $X$ and $Y$ be random vectors of dimension
$p$ and $q$. If there is a
matrix $\beta\in\real^{p \times d}$ with $d < p$ such that
%
\begin{equation}
\label{eqsdr} Y \indep X | \beta^{\mathsf{T}}X,
\end{equation}
then the subspace spanned by the columns of $\beta$ is called a
\textit{sufficient dimension reduction} (\textit{SDR}) \textit{subspace}. Under mild
conditions, the intersection of all such subspaces still satisfies
(\ref{eqsdr}), and is called the \textit{central subspace}, denoted by
$\lcentral$; see Li (\citeyear{Li91N2}, \citeyear{Li92}), \citet{LiDua89},
\citet{DuaLi91}, \citet{Li91N1}, Cook (\citeyear{Coo94},
\citeyear{Coo98N1}). A general condition for the existence of the central subspace is
given by \citet{YinLiCoo08}.

Several recent papers have combined sufficient dimension reduction and
kernels; see \citet{Aka01}, \citet{BacJor03},
\citet{FukBacGre07}, \citet{Wu08}, \citet{WuLiaMuk08},
\citet{HsiRen09}, \citet{YehHuaLee09}, \citet{ZhuLi11}
and \citet{LiArtLi}. This proliferation of work, in
addition to producing versatile methods for extracting nonlinear
sufficient predictors, points toward a general synthesis between the
notions of sufficiency at the core of SDR and the ability to encompass
nonlinearity afforded by kernel mappings. To achieve this synthesis,
explore its many ramifications and broad scope and develop new
estimators based on it, are the goals of this paper.

Specifically, we articulate a general formulation that comprises both
linear and nonlinear SDR, and parallels the basic theoretical
developments pioneered by Li (\citeyear{Li91N2}, \citeyear{Li92}) and
Cook (\citeyear{Coo94}, \citeyear{Coo98N2}, \citeyear{Coo98N1}). This
formulation allows us to study linear and nonlinear SDR comparatively
and, somewhat surprisingly, to \textit{relax} some stringent conditions
required by linear SDR. For example, a linear conditional mean
[\citet{Li91N2}, \citet{Coo98N1}] is no longer needed for
unbiasedness, and the sufficient conditions for existence and
uniqueness of the central subspace are far more general and
transparent. Finally, our general formulation links linear and
nonlinear SDR to the classical notions sufficiency, completeness and
minimal sufficiency, which brings insights and great clarity to the SDR
theory.

Our developments and the sections of this paper, can be summarized as
follows. In Section~\ref{sec2}, we build upon
the ideas of \citet{Coo07} and \citet{LiArtLi} to define an
SDR $\sigma$-field as a sub $\sigma$-field $\G$ of $\sigma(X)$
(the $\sigma$-field generated by $X$) such that
$
Y \indep X | \G,
$
and the corresponding SDR class as the set of all square-integrable,
$\G
$-measurable functions.
Under very mild conditions---much milder than the corresponding
conditions for linear SDR [\citet{YinLiCoo08}]---we show that
there exists a unique minimal $\sigma$-field $\cfield$ that satisfies
$Y \indep X | \cfield$, which we call the \textit{central $\sigma$-field}.
The set of all $\cfield$-measurable,
square-integrable functions is named the \textit{central
class}.\vadjust{\goodbreak}

In Section~\ref{sec3}, we provide two additional definitions that generalize
concepts in Cook (\citeyear{Coo98N1}), \citet{LiZhaChi05} and
\citet{LiArtLi}: a~class of functions is \textit{unbiased} if its
members are $\cfield$-measurable, and \textit{exhaustive} if they generate
$\cfield$. Next, we show that the special class
%
\begin{equation}
\label{eqcore} {L_2 (P_{ X})}\ominus\bigl[ {L_2
(P_{ X})}\ominus{L_2 (P_{ Y})}\bigr]
\end{equation}
is unbiased, where ${L_2 (P_{ X})}$ and ${L_2 (P_{ Y})}$ are the spaces of
square-integrable functions of $X$ and $Y$. For reasons detailed in
Section~\ref{sec3}, we call
this class the \textit{regression class}.

In Section~\ref{sec4}, we introduce the \textit{complete dimension reduction
class}: if
$\gfield\subseteq\sigma(X)$ is a $\sigma$-field and for each $\G
$-measurable $f \in{L_2(P_{X})}$ we have
\[
E\bigl(f(X) | Y\bigr) = 0 \qquad\mbox{almost surely} \quad\implies\quad f(X) = 0
\qquad\mbox{almost surely},
\]
then we say that the class of $\gfield$-measurable functions in
${L_2(P_{X})} $ is complete. We prove that when a complete sufficient
dimension reduction (CSDR) class exists it is unique and coincides with
the central class. We further show that the CSDR class coincides with
the regression class---which is therefore not just unbiased, but also
exhaustive.

In Section~\ref{sec5} we establish a critical relationship between the
regression class and a covariance operator
linking $X$ and $Y$ and, based on this, we generalize sliced inverse
regression [SIR; \citet{Li91N2}]
to a method (GSIR) that can recover the regression class---and hence
is unbiased and exhaustive under completeness.
In Section~\ref{sec6}, we consider the case where the central class is not complete,
so that GSIR is unbiased but no longer exhaustive. By introducing a
novel conditional variance operator,
we generalize sliced average variance estimation [SAVE, \citet{Li91N1}] to a method (GSAVE) that
can recover a class larger
than the regression class. 
Here, the situation is similar to that in the linear SDR setting, where
it is well known that
%
\begin{equation}
\label{eqcookcritchley} \mbox{SIR subspace} \subseteq\mbox{SAVE subspace}
\subseteq{\cal S}_{ Y|X};
\end{equation}
see Cook and Critchley (\citeyear{CoCr00}),
Ye and Weiss (\citeyear{YeWe03}), \citet{LiZhaChi05} and \citet{LiWan07}.

In Section~\ref{sec7} we develop algorithms for the sample versions of GSIR and
GSAVE, and a cross-validation algorithm to determine
regularizing parameters. In Section~\ref{sec8} we compare GSIR and GSAVE with
some existing methods by simulation and on actual data sets.
Section~\ref{sec9} contains some concluding remarks. Some highly technical
developments are provided in the supplementary material [\citet{LeeLiChi}].

\section{\texorpdfstring{Sufficient dimension reduction $\sigma$-fields and classes}{Sufficient dimension reduction sigma-fields and classes}}\label{sec2}

Let $(\om, \field, P)$ be a probability space and $(\omx, \fieldx)$,
$(\omy, \fieldy)$
and $(\omxy, \fieldxy)$ be measurable spaces. For convenience, assume
that $\omxy= \omx\times\omy$ and
$\fieldxy= \fieldx\times\fieldy$. Let $X$, $Y$ and $(X,Y)$ be random\vadjust{\goodbreak}
elements that take values
in $\omx$, $\omy$ and $\omxy$, with distributions
$
\px$, $\py$, $\pxy$, which are dominated by $\sigma$-finite
measures. Let
\[
\sigma(X) = X^{-1}(\borelx),\qquad \sigma(Y) = Y^{-1}(\borely),\qquad
\sigma(X,Y) = (X,Y)^{-1}(\fieldxy),
\]
and finally let $\pygx(\cdot|\cdot)\dvtx
\borely\times\ox\to\R$
be the conditional distribution of $Y$ given~$X$.


\begin{definition}
\label{definitionsufficiency}
A sub $\sigma $-field $\G$ of $\sigma(X)$ is an SDR $\sigma$-field for
$Y$ versus $X$ if it satisfies
%
\begin{equation}
\label{eqsufficientsigmafield} Y \indep X | \G,
\end{equation}
that is, if $Y$ and $X$ are independent given $\G$.
\end{definition}

This definition is sufficiently general to accommodate the two cases of
nonlinear sufficient dimension reduction that
interest us the most. The first case is when $\omx= \real^p$ and
$\omy
= \real^q$ for
some positive integers $p$ and $q$, and $\fieldx$, $\fieldy$ and
$\fieldxy$ are Borel $\sigma$-fields
generated by the open sets in $\real^p$, $\real^q$ and $\real^{p+q}$.
Clearly, in this case,
the conditional independence (\ref{eqsufficientsigmafield}) is a
generalization of
(\ref{eqsdr}) for linear SDR: if we take $\G= \sigma(\beta^{\mathsf{T}}
X)$, then (\ref{eqsufficientsigmafield})
reduces to (\ref{eqsdr}). 


The second case is when $X$ or $Y$, or both of them, are random
functions. In this case
Definition~\ref{definitionsufficiency} is a generalization of the
linear SDR
for functional data introduced by Ferr\'{e} and Yao (\citeyear{FeYa03}),
and \citet{HsiRen09}. Specifically,
let $[a,b]$ be a closed interval, $\lambda$ the Lebesgue measure and
$L_2 (\lambda)$ the class of functions on $[a,b]$ that are square
integrable with
respect to $\lambda$. Let $\omx= L_2 (\lambda)$ and $\omy= \real$. In
this case, each
$X(\omega)$ is a function in $L_2 (\lambda)$, which, depending on
applications,
could be, say, a growth curve or the fluctuation of a stock price.
Let $h_1,\ldots, h_d$ be functions in $L_2 (\lambda)$.
Ferr\'{e} and Yao (\citeyear{FeYa03}) considered
the following functional dimension reduction problem:
%
\begin{equation}
\label{eqferreyao1} Y \indep X | \langle X, h_1
\rangle_{ L_2 (\lambda)},\ldots, \langle X, h_d \rangle_{ L_2 (\lambda)}.
\end{equation}
This generalizes linear SDR
to the infinite-dimensional case, but not to the nonlinear case, because
$\langle X, h_1 \rangle_{ L_2 (\lambda)},\ldots, \langle X, h_d
\rangle_{ L_2 (\lambda)}$ are linear in $X$.
\citet{HsiRen09} considered a more general setting
where the sample paths $\{X_t(\omega)\dvtx  t \in J\}$ need not lie within
$L_2 (\lambda)$. Still, their generalization is inherently linear in
the same sense that problem (\ref{eqferreyao1}) is linear. In
contrast, our formulation in (\ref{eqsufficientsigmafield})
allows
an arbitrary sub $\sigma$-field of $\sigma(X)$, which need not be
generated by linear
functionals. Interestingly, as we will see Section~\ref{sec5}, it is the
relaxation of linearity that allows us to remove
a restrictive linear conditional mean assumption 
used both in Ferr\'{e} and Yao (\citeyear{FeYa03}) (Theorem~2.1), and in \citet{HsiRen09}, assumption (IR2).


The notion of sufficiency underlying SDR, as defined by (\ref{eqsdr})
and (\ref{eqsufficientsigmafield}), is different from the classical
notion of sufficiency because $\cal G$ is allowed to depend on any
parameter in the joint distribution of $\pxy$. For example, $\G=
\sigma(\beta^{\mathsf{T}}X)$ depends on the parameter $\beta$ [or
rather, the meta-parameter $\sp(\beta)$] which characterizes the
conditional distribution of $Y|X$. Nevertheless, both notions imply a
reduction, or simplification, in the representation of a stochastic
mechanism---the SDR one through a newly constructed predictor, and the
classical one through a statistic. Indeed, it is partly by exploring
and exploiting this similarity that we developed our theory of
nonlinear SDR.


Obviously there are many sub $\sigma$-fields of $X$ that satisfy (\ref
{eqsufficientsigmafield}),
starting with $\sigma(X)$ itself---which induces no reduction.
For maximal dimension reduction we seek the smallest
such $\sigma$-field. As in the case of classical sufficiency, the
minimal $\mathrm{SDR}$ \mbox{$\sigma$-field} does not universally exist, but exists
under very mild assumptions. The next theorem gives the sufficient condition
for the minimal $\mathrm{SDR}$ $\sigma$-field to uniquely exist. The
proof echoes
\citet{Bah54}, which established
the existence of the minimal sufficient $\sigma$-field in the
classical setting.

\begin{theorem}\label{theoremexistenceanduniqueness}
Suppose that the family of probability measures $\{
P_{ X | Y} (\cdot|y)\dvtx\break  y \in
\Omega_{ Y} \}$ is dominated by a $\sigma$-finite measure. Then
there is a unique sub $\sigma$-field
$\G^*$ of $\sigma(X)$ such that:
\begin{longlist}[(2)]
\item[(1)] $Y \indep X | \G^*$;
\item[(2)] if $\G$ is a sub $\sigma$-field of $\sigma(X)$ such that $Y
\indep X | \G$, then $\G^* \subseteq\G$.
\end{longlist}
\end{theorem}


\begin{pf}
Let $\Pi_y=P_{ X|Y} (\cdot|y)$ and $\pmathbb= \{\Pi_y\dvtx y \in
\omy
\}$.
Since $\mathbb P$ is dominated by a $\sigma$-finite measure,
it contains a countable subset
$
\qmathbb= \{Q_k\dvtx  k=1, 2, \ldots\}
$ such that $\qmathbb\equiv\pmathbb$, where $\equiv$ means two
families of measures dominating
each other.
Let $\{c_k\dvtx  k = 1,2,\ldots\}$ be a sequence of positive numbers that
sum to 1,
and let $Q_{0} = \sum_{k=1}^\infty c_k Q_k$. Then $Q_{0}$
is a probability measure on $\Omega_{ X}$
such that
$\{Q_{0}\} \equiv\qmathbb\equiv\pmathbb$.
Let $\pi_y = d \Pi_y / d Q_{0}$ and $\G$ be a sub $\sigma$-field
of $\sigma(X)$.
We claim that the following statements are equivalent:
\begin{longlist}[(2)]
\item[(1)] $Y \indep X | \G$;
\item[(2)] $\pi_y$ is essentially measurable with respect to $\G$ for
all $y \in\Omega_{ Y}$ modulo~$Q_{0}$.
\end{longlist}

\textsc{Proof of 1 $\Rightarrow$ 2}.\quad Let $B \in\fieldx$. Then
\begin{eqnarray*}
E_{ Q_0} \bigl(\pi_y (X) I_{ B} (X) \bigr) &=&
E_{\Pi_y} \bigl( I_{
B} (X) \bigr) = E_{\Pi_y} \bigl[
E_{\Pi_y} \bigl( I_{ B} (X) | \G\bigr) \bigr] \\
&=& E_{ Q_0}
\bigl[ E_{\Pi_y} \bigl( I_{ B} (X) | \G\bigr)
\pi_y(X) \bigr].
\end{eqnarray*}
By 1,
$\Pi_y (B| \G)$ is the same for all $y \in\Omega_{ Y}$. Hence
$\Pi_y (B| \G) = Q_k (B| \G)$ for all $k$,
which implies
$\Pi_y (B| \G) = Q_{0} (B| \G)$. Hence we can rewrite the
right-hand side of the above equalities as
\[
E_{ Q_0} \bigl[ E_{ Q_0} \bigl( I_{ B} (X) | \G\bigr)
\pi_y(X) \bigr] = E_{ Q_0} \bigl[ I_{ B} (X)
E_{ Q_0} \bigl(\pi_y(X)| \G\bigr)\bigr].
\]
Thus the following equality holds for all $B \in\fieldx$:
\[
E_{ Q_0} \bigl(\pi_y (X) I_{ B} (X) \bigr) =
E_{ Q_0} \bigl[ I_{
B} (X) E_{ Q_0} \bigl(
\pi_y(X)| \G\bigr)\bigr],
\]
which implies $\pi_y (X)= E_{ Q_0} (\pi_y(X)| \G)$ a.s.
$Q_{0}$.

\textsc{Proof of 2 $\Rightarrow$ 1}.\quad 
For any $A \in{\cal G}$,
%
\begin{eqnarray*}
E_{\Pi_y} \bigl[E_{Q_0}\bigl( I_{ B}(X) |\G\bigr)
I_{ A} (X) \bigr] &=& E_{ Q_0} \bigl[ E_{ Q_0}\bigl(
I_{ B}(X) |\G\bigr) I_{ A} (X) \pi_y(X) \bigr]
\\
&=& E_{ Q_0} \bigl[ I_{ B}(X)I_{ A} (X)
E_{ Q_0}\bigl( \pi_y(X) |\G\bigr) \bigr].
\end{eqnarray*}
By 2, $E_{Q_{0}} ( \pi_y (X) | \G) = \pi_y (X)$. Hence the right-hand
side becomes
\[
E_{ Q_0} \bigl[ I_{ B}(X)I_{ A} (X)
\pi_y(X) \bigr] = E_{
\Pi
_y} \bigl[ I_{ B}(X)I_{ A}
(X) \bigr] = \Pi_y ( X \in A \cap B ).
\]
Thus $E_{ Q_0}(I_{ B}(X) | \G) = Q_{0}(B | \G)$ is the
conditional probability $\Pi_y(B|\G)$,
which means $\Pi_y(B|\G)$ does not depend on $y$.
That is, 1 holds.

Now let $\G^*$ be the intersection of all $\mathrm{SDR}$ $\sigma
$-fields $\G$.
Then $\G^*$ is itself
a $\sigma$-field. Moreover, since $\pi_y$ is essentially measurable
with respect to all $\mathrm{SDR}$ $\sigma$-fields
for all $y \in\Omega_{ Y}$, it is also essentially measurable
with respect to $\G^*$ for all $y \in\Omega_{ Y}$.
Consequently, $\G^*$ is itself an $\mathrm{SDR}$ $\sigma$-field, which implies
that it is also the smallest $\mathrm{SDR}$ $\sigma$-field.
If $\G^{**}$ is another smallest $\mathrm{SDR}$ $\sigma$-field, then we know
$\G^* \subseteq\G^{**}$ and $\G^{**} \subseteq\G^*$. Thus $\G^*$ is
unique.
\end{pf}

We can now naturally introduce the following definition:

\begin{definition}
Suppose that the class of probability measures $\{P_{
X|Y}(\cdot|y)\dvtx \break y \in\Omega_{ Y} \}$ on $\Omega_{ X}$ is dominated
by a $\sigma$-finite measure. Then we call the $\sigma$-field $\G^*$ in
Theorem \ref {theoremexistenceanduniqueness} the central $\sigma$-field
for $Y$ versus $X$, and denote it by $\G_{ Y|X}$.
\end{definition}


Notably, this set up characterizes dimension reduction solely in terms
of conditional independence.
However, explicitly turning to functions and introducing an additional
mild assumption of square
integrability are very consequential for further
development because they allow us to work with structures such as
orthogonality and projection.


Let $L_2 (\pxy)$, $L_2 (\px)$ and $L_2 (\py)$ be the spaces of functions
defined on $\Omega_{ XY}$,
$\Omega_{ X}$ and $\Omega_{ Y}$ that are square-integrable
with respect to $P_{ XY}$,
$ P_{ X}$ and $P_{ Y}$, respectively. Since constants are
irrelevant for dimension reduction,
we assume throughout that all functions in ${L_2(P_{X})}$,
${L_2(P_{Y})}$ and
$L_2 (P_{ XY})$ have mean 0.
Given a sub $\sigma$-field $\G$ of $\sigma(X, Y)$, we use $\M_\G$ to
denote the class of all functions
$f$ in $L_2 (P_{ XY})$ such that $f(X)$ is $\G$-measurable. If
$\G
$ is generated by a random vector,
say $X$, then we use $\M_{ X}$ to abbreviate $\M_{\sigma{(X)}}$.
It can be easily shown that, for any $\G$, $\M_\G$ is a linear subspace
of $L_2 (P_{ XY})$.

\begin{definition}
Let $\G$ be an $\mathrm{SDR}$ $\sigma$-field and $\G_{ Y|X}$ be the
central $\sigma$-field. Then $\M_\G$ is called an $\mathrm{SDR}$ class,
and $\M_{\G _{Y|X}}$ is called the central class. The latter class is
denoted by $\S_{ Y|X}$.
\end{definition}


The central class, comprising square-integrable functions that are
measurable with respect to the
central $\sigma$-field $\G_{ Y|X}$, represents our generalization
of the central space
$\spc_{ Y|X}$ defined in linear SDR; see the \hyperref[sec1]{Introduction}.

\section{Unbiasedness and exhaustiveness}\label{sec3}\label{sectionunbiasedcharacterization}




In linear SDR, the goal is to find a set of vectors that span $\spc_{ Y|X}$. If a matrix $\gamma$ satisfies $\sp(\gamma) \subseteq
\spc_{ Y|X}$, we say that $\gamma$ is unbiased [\citet{Coo98N1}].
If $\sp(\gamma) = \spc_{ Y|X}$, we say that $\gamma$ is exhaustive
[\citet{LiZhaChi05}].
Note that when $\sp(\gamma) \subseteq\spc_{ Y|X}$, $\gamma^{\mathsf{T}}
X$ is a linear function of
$\beta^{\mathsf{T}}X$, where $\beta$ is any matrix such that $\sp
(\beta) =
\spc_{ Y|X}$; if $\sp(\gamma)
= \spc_{ Y|X}$, then $\gamma^{\mathsf{T}}X$ is an injective linear
transformation of $\beta^{\mathsf{T}}X$.
In the nonlinear setting, we follow the same logic but remove the
linear requirement. Part of the following definition was given in
\citet{LiArtLi}.

\begin{definition}
A class of functions in ${L_2(P_{X})}$ is unbiased for $\central$ if
its members are $\cfield$-measurable, and exhaustive for $\central$ if
its members generate $\cfield$.
\end{definition}

Next, we look into what type of functions are unbiased. The lemma below
provides a characterization of
the orthogonal complement of $\M_{\G}$ that will be used many times in
the subsequent development.
Its proof is essentially the definition of the conditional expectation,
and is omitted.

\begin{lemma}\label{lemmakeylemma} Suppose $U$ is a random element
defined on $(\sample,\field)$,
$\G$ is a sub $\sigma$-field of $\sigma(U)$ and $f \in L_2 (P_{
U})$. Then $f$ is orthogonal
to $\M_\G$ ($f \perp\M_\G$) if and only if $E[ f (U) |\G] = 0$.
\end{lemma}

%

Note that $\indep$ and $\perp$ have different meanings: the former
means independence; the latter means orthogonality. For two subspaces,
say $\spc_1$ and $\spc_2$, of a generic Hilbert space $\cal H$, we use
$\spc_1 \ominus\spc_2$ to denote the subspace $\spc_1 \cap\spc_2^\perp
$. The following theorem explicitly specifies a class of functions,
which we call \textit{regression class}, that is unbiased for
$\central$.

\begin{theorem}\label{theoremsufficiencyalone} If the family $\{\Pi_y
\dvtx  y \in\Omega_{ Y} \}$
is dominated by a $\sigma$-finite measure, then
%
\begin{equation}
\label{equnbiasedness} L_2 (P_{ X}) \ominus\bigl[
L_2 (P_{ X}) \ominus L_2 (P_{
Y})
\bigr] \subseteq\central.
\end{equation}
\end{theorem}
\begin{pf}
It is equivalent to show that $
L_2 (P_{ X}) \ominus\central\subseteq
L_2 (P_{ X}) \ominus L_2 (P_{ Y}).
$
If $f \in L_2 (P_{ X})\ominus\central$, then, by Lemma \ref
{lemmakeylemma},
$E [f (X) | \cfield] = 0$.
Since $\cfield$ is a sufficient $\sigma$-field,
\[
E \bigl[ f (X) | Y \bigr] = E \bigl[E \bigl(f (X) | Y, \cfield\bigr) | Y \bigr]
= E \bigl[E \bigl(f (X) |\cfield\bigr) | Y \bigr] =0.
\]
By Lemma~\ref{lemmakeylemma} again, $f \perp\M_{ Y} $. Because
$\M_{ Y} =L_2 (P_{ Y})$,
we have $f \in L_2 (P_{ X}) \ominus L_2 (P_{ Y})$.\vadjust{\goodbreak}
\end{pf}

The intuition behind the term ``regression class'' is that ${L_2(P_{X})}
\ominus{L_2(P_{Y})}$ resembles the residual
in a regression problem; thus ${L_2(P_{X})}\ominus[
{L_2(P_{X})}\ominus{L_2(P_{Y})}
]$ is simply the orthogonal complement
of the ``residual class.'' Henceforth we write the regression class as
$\complete$.

\section{Complete and sufficient dimension reduction classes}\label{sec4}

After showing that the regression class (\ref{eqcore}) is unbiased, we
investigate under what conditions it is also exhaustive for the central class
$\central$. To this end we need to introduce the notion of complete
classes of functions in ${L_2(P_{X})}$.

\begin{definition}\label{definitioncomplete}
Let $\G\subseteq \sigma(X)$ be a sub $\sigma$-field. The class $\M_\G$
is said to be complete if, for any $g \in\M_\G$,
\[
E \bigl[g (X) | Y \bigr] = 0 \qquad\mbox{a.s. } P \quad\Rightarrow\quad
g (X) = 0 \qquad\mbox{a.s. } P.
\]
\end{definition}


Again there are similarities and differences between completeness as
defined here and in the classical
setting. A complete and sufficient statistic in the classical setting
is a rather restrictive
concept, often associated with exponential families, the uniform
distribution, or the order statistics.
In contrast, completeness here is a rather general concept. To
demonstrate this point, in the next two propositions we give two
examples of complete and sufficient dimension reduction classes. In particular,
the first shows that if $Y$ is related to $X$ through \textit{any}
regression model with additive error, then the subspace of
$L_2 (P_{ X})$ determined by the regression function is a complete
and sufficient dimension reduction class. In the following, $[{L_2(P_{X})}
]^q$ denotes the $q$-fold Cartesian product of ${L_2(P_{X})}$.

\begin{proposition}\label{prop1} Suppose there exists a function $h \in[{L_2(P_{X})}]^q$
such that
%
\begin{equation}
\label{eqgeneralregression} Y = h(X) + \varepsilon,
\end{equation}
where $\varepsilon\indep X$ and $E(\varepsilon) = 0$. Then $\M_{h(X)}$
is a complete and sufficient
dimension reduction class for $Y$ versus $X$.
\end{proposition}

Note that, since ${L_2(P_{X})}$ is centered, we have implicitly assumed
that
\mbox{$E[h(X)] = 0$} [and hence
$E(Y) = 0$]. However, this does not entail any real loss of generality
because the proof below
can be easily modified for the case where ${L_2(P_{X})}$ is not centered.

\begin{pf*}{Proof of Proposition~\ref{prop1}} Suppose $m \in\M_{h(X)} $ and $E[ m (X) | Y]
= 0 $ a.s.~$P$.
Then there is a measurable
function $g\dvtx  \real^q \to\real$
such that $m = g \circ h $. Let \mbox{$U = h(X)$}. Then $E( g (U) | Y) = 0 $
a.s. $P$.
By Lemma~\ref{lemmakeylemma}, for any $f \in L_2 (P_{ Y})$, we have
$E[ g(U) f(Y) ] = 0$.
In particular,
$
E[ g(U) e^{it^{\mathsf{T}}Y} ] = 0.
$
Because $U \indep\varepsilon$, this
implies
\[
E\bigl[ g(U) e^{i t^{\mathsf{T}}U} \bigr] E \bigl(e^{i t^{\mathsf{T}}\varepsilon}\bigr) = E\bigl[
g(U) e^{i t^{\mathsf{T}}U} e^{i t^{\mathsf{T}}\varepsilon} \bigr] = E\bigl[ g(U) e^{i t^{\mathsf{T}}Y}
\bigr]=0.
\]
Hence $E[ g(U) e^{i t^{\mathsf{T}}U} ] = 0$. By the uniqueness of
inverse Fourier transformation we see that $g(U) = 0$ a.s. $P$, which
implies $m(X) = (g \circ h) (X) = 0 $ a.s. $P$.~%
\end{pf*}

The expression in (\ref{eqgeneralregression}) covers many useful
models in statistics and
econometrics. For example, any homoscedastic parametric or
nonparametric regression, such as
the single index and the multiple index models [\citet{IchLee91},
\citet{HarHalIch93}, \citet{YinLiCoo08}], are special cases of (\ref
{eqgeneralregression}).
Thus, complete and sufficient dimension reduction classes exist for all
those settings. The next
proposition considers a type of inverse regression model, in which $X$
is transformed into two
components, one of which is related to $Y$ by an inverse linear
regression model, and the other
independent of the rest of the data.

\begin{proposition} Suppose $q<p$, $\Omega_{ Y}$ has a nonempty
interior, and $P_{ Y}$
is dominated by the Lebesgue measure on $\real^q$. Suppose there exist
functions $g \in[{L_2(P_{X})}]^q$
and $h \in[{L_2(P_{X})}]^{p-q}$ such that:
\begin{longlist}[(4)]
\item[(1)] $g(X) = Y + \varepsilon$, where $Y \indep\varepsilon$, and
$\varepsilon\sim N(0, \Sigma)$;
\item[(2)] $\sigma(g(X), h(X)) = \sigma(X)$;
\item[(3)] $h(X) \indep(Y, g(X))$;
\item[(4)] the induced measure $\px\circ g^{-1}$ is dominated by the
Lebesgue measure on~$\real^q$.
\end{longlist}
Then $\M_{g(X)} $ is a complete sufficient dimension reduction class
for $Y$ versus~$X$.
\end{proposition}


\begin{pf}
Assumption 3 implies $Y \indep h(X) | g(X)$, which, by
assumption~2, implies $Y \indep X | g(X)$.
That is, $\M_{g(X)} $ is an SDR class. Let $u\in\M_{g(X)}$. Then $u
= v
\circ g$
for some measurable function $v\dvtx  \real^q \to\real$. Let $U = g(X)$.
Suppose that $E[v (U) | Y] = 0 $
almost surely $P$. 
Because $Y \indep\varepsilon$, this implies $\py(\{y\dvtx  E v(y +
\varepsilon) = 0 \} ) = 1$.
In other words,
\[
\int_{\real^q} v(t)\frac{1}{(2 \pi)^{q/2}| \Sigma|^{1/2}} e^{-
(t-y)^{\mathsf{T}}\Sigma^{-1}(t-y) / 2} \,d t = 0
\]
a.s. $\py$.
This implies
\[
\int v(t) e^{-t^{\mathsf{T}}\Sigma^{-1}t/2} e^{y^{\mathsf{T}}\Sigma
^{-1}t} \,d t = 0 \quad\implies\quad \int v(\Sigma s )
e^{-s^{\mathsf{T}}\Sigma s/2} e^{y^{\mathsf{T}}s} \,d s = 0
\]
a.s. $\py$,
where $s = \Sigma^{-1}t$.
Because $\Omega_{ Y}$ contains an open set in $\real^q$ and the
above function of $y$ is analytic,
by the analytic continuation theorem, the above function is 0
everywhere on $\real^q$.
Hence, by the uniqueness of inverse Laplace transformation, we have
\[
v(\Sigma s ) e^{-s^{\mathsf{T}}\Sigma s/2} = 0 \qquad\mbox{almost surely $\lambda$},
\]
where $\lambda$ is the Lebesgue measure on $\real^q$. But, because
$e^{-s^{\mathsf{T}}\Sigma s/2} > 0$, we have
$v(\Sigma s ) = 0$ a.s. $\lambda$ or equivalently $v(t) = 0$ a.s.\vadjust{\goodbreak}
$\lambda$. By the change of variable theorem,
\[
\int_{v\circ g (x) \ne0} d \px= \int_{v(t) \ne0} d \px
\circ g^{-1}.
\]
By assumption 4, \mbox{$\px\circ g^{-1}\ll\lambda$}. Hence the above
integral is 0, implying $v\circ g(x) = 0$ a.s. $\px$, or, equivalently,
$v\circ g(X) = 0$ a.s. $P$.
\end{pf}

Inverse regressions of this type are considered in \citet{Coo07},
\citet{CooFor09}, and \citet{CooLiChi10} for linear SDR. The
above two propositions show that a complete and sufficient dimension
reduction class exists for a reasonably wide range of problems,
including forward and inverse regressions of very general,
nonparameterized form. The next theorem shows that when a complete and
sufficient dimension reduction class exists, it is unique and coincides
with the central class. Once again, the situation here echoes that in
classical theory, where a complete and sufficient statistic, if it
exists, coincides with the minimal sufficient statistic; see
\citet{Leh81}.

\begin{theorem} Suppose $\{\Pi_y\dvtx  y \in\Omega_{ Y}\}$ is
dominated by a
$\sigma$-finite measure, and $\G$ is a sub $\sigma$-field of $\sigma(X)$.
If $\M_\G$ is a complete and sufficient dimension reduction class, then
\[
\M_\G= \complete= \S_{ Y|X}.
\]
\end{theorem}

\begin{pf}
If $f \perp\complete$, then by Lemma
\ref{lemmakeylemma}, $E(f | Y) = 0$ which, because $\M_\G$ is
sufficient, implies
\[
E\bigl[ E( f | \G) | Y \bigr] = 0.
\]
Because $\M_\G$ is complete and because $E(f | \G) \in\M_\G$, we have
$E(f | \G) = 0$. By Lemma
\ref{lemmakeylemma}, this implies $f \perp\M_\G$. Thus we have
proved $\M_\G\subseteq\complete$.
However, by Theorem~\ref{theoremsufficiencyalone} we know that
$\complete\subseteq\central\subseteq\M_\G$.
This proves the desired equality.
\end{pf}

%

%


\section{Generalizations of SIR and their population-level
properties}\label{sec5}\label{sectiongsir}

From the previous developments we see that the subspace
$
L_2 (P_{ X}) \ominus L_2 (P_{ Y})
$
of $L_2 (P_{ X})$ plays a critical role in nonlinear SDR. Its
orthogonal complement in
$L_2 (P_{ X})$ coincides with the central class $\central$ under
completeness, and even without completeness
it is guaranteed to be inside $\central$. It turns out that this
subspace can be expressed as the range of a certain bounded linear
operators. This representation ensures that
estimation procedures can rely on simple spectral decompositions,
rather than complicated
numerical optimizations. We first introduce some covariance operators
which are the building block
of this approach.

\subsection{Covariance operators}\label{sec5.1}

Since constants are irrelevant here (e.g., $f$ and $f + 3$ can be
considered as the\vadjust{\goodbreak}
same function), we will speak of set relations modulo constants. If $A$
and $B$ are sets,
then we say $A \subseteq B$ modulo constants if for each $f \in A$
there is $c \in\real$ such that
$f + c \in B$. We say that $A$ is a dense subset of $B$ modulo
constants if (i) $A \subseteq B$ modulo constants and (ii) for each
$f \in B$, there is a sequence $\{f_n\} \subseteq A$ and a sequence of
constants $\{c_n \}
\subseteq\real$ such that $\{f_n + c_n \} \subseteq A$ and $f_n + c_n
\to f$ in the topology
for $B$. Let $\hx$ and $\hy$ be Hilbert spaces of functions of $X$ and
$Y$ satisfying the
conditions:
\begin{longlist}[(B)]
\item[(A)] $\hxquo$ and $\hyquo$ are dense subsets of ${L_2 (P_{
X})}$ and
${L_2 (P_{ Y})}$ modulo constants;
\item[(B)] there are constants $C_1 > 0$ and $C_2 > 0$ such that $\var
[f(X)] \le C_1 \|f \|_\hx$ and
$\var[g(Y)] \le C_2 \| g \|_\hy$.
\end{longlist}
%
Although we will later take $\hx$ and $\hy$ to be reproducing kernel
Hilbert spaces (RKHS), our theory is not restricted to such spaces. In
particular, we do not
require the evaluation functionals
[such as $f \mapsto f(x)$ from $\hx$ to~$\real$] to be continuous.

For two generic Hilbert spaces $\hone$ and $\htwo$, let
${\cal B}(\hone,\htwo)$ denote the class of all bounded linear
operators from
$\hone$ to $\htwo$, and let
${\cal B}(\hone)$ abbreviate ${\cal B}(\hone,\hone)$. We denote the
range of a
linear operator $A$ by $\ran A$, the kernel of $A$
by $\ker A$, and the closure of $\ran A$ by $\overran A$.
Under assumption (B), the symmetric bilinear form $u\dvtx  \hx\times\hx
\to\real$ defined by
$u(f,g) = \cov[ f(X), g(X) ]$ is bounded and thus induces
an operator $\mxx\in{\cal B}( \hx)$ that satisfies
$
\langle f, \mxx g \rangle_\hx= u(f,g).
$ Similarly, the bounded bilinear form $(f,g) \mapsto\cov[f(Y), g(Y)]$
from $\hy\times\hy$ to $\real$
defines an operator $\myy\in{\cal B}(\hy)$.
Let ${\cal G}_{X}$ and ${{\cal G}_{Y}}$ represent the subspaces
$\overran\mxx$
and $\overran\myy$.

%
%


%
%
%


\begin{definition} Suppose conditions (A) and (B) are satisfied. We
define the covariance operators
$\coxx\dvtx  {\cal G}_{X}\to{\cal G}_{X}$, $\coyy\dvtx  {{\cal G}_{Y}}\to
{{\cal G}_{Y}}$ and $\coyx\dvtx  {\cal G}_{X}\to{{\cal G}_{Y}}$
through the relations
\begin{eqnarray*}
\langle f, \coxx g \rangle_{{\cal G}_{X}}&=& \langle f, g \rangle_{{L_2 (P_{
X})}},\qquad
\langle f, \coyy g \rangle_{{\cal G}_{Y}}= \langle f, g \rangle_{{L_2 (P_{
Y})}},\\
\langle f, \coyx g \rangle_{{\cal G}_{Y}}&=& \langle f, g \rangle_{{L_2 (P_{ Y})}}.
\end{eqnarray*}
\end{definition}

These operators are essentially the same as those introduced by
Fukumizu, Bach and Jordan (\citeyear{FukBacJor03}, \citeyear{FukBacJor09}),
except that here we do not assume $\hx$ and $\hy$ to be RKHS.
By Baker [(\citeyear{Bak73}), Theorem 1], there is a unique operator $\ryx\in\B(
{{\cal G}_{X}}, {{\cal G}_{Y}})$ such that
$\coyx= \coyy^{1/2}\ryx\coxx^{1/2}$. We call $\rxy$ the \textit{correlation operator}.
In order to connect these operators with the central class, which is an
${L_2(P_{X})}$-object, we need to
extend the domains of
$\coxx^{1/2}$ and $\coyy^{1/2}$ from ${{\cal G}_{X}}$ and ${{\cal G}_{Y}}$ to
${L_2(P_{X})}$ and
${L_2(P_{Y})}$.
The following extension theorem is important and nontrivial, but since
the material presented here can be
understood without its proof we relegate it to
the supplementary material [\citet{LeeLiChi}].\vadjust{\goodbreak}

\begin{theorem} Under assumptions \textup{(A)} and \textup{(B)}, there exist unique isomorphisms
\[
\tilde\Sigma_{ XX}^{1/2}\dvtx {L_2(P_{X})}
\to{\cal G}_{X}, \qquad\tilde\Sigma_{ YY}^{1/2}
\dvtx {L_2(P_{Y})}\to{{\cal G}_{Y}}
\]
that agree with $\Sigma_{ XX}^{1/2}$ and $\Sigma_{ YY}^{1/2}$ on
${\cal G}_{X}$ and ${{\cal G}_{Y}}$ in the sense
that, for all $f \in{\cal G}_{X}$ and $g \in{{\cal G}_{Y}}$,
\[
\tilde\Sigma_{ XX}^{1/2}(f - E f) = \Sigma_{ XX}^{1/2}f,\qquad
\tilde\Sigma_{ YY}^{1/2}( g - E g) = \Sigma_{ YY}^{1/2}g.
\]
Furthermore, for any $f \in{L_2(P_{X})}$, $g \in{L_2(P_{Y})}$ we have
%
\begin{equation}
\label{eqtarget1} \bigl\langle\tilde\Sigma_{ YY}^{1/2}g, \ryx
\tilde\Sigma_{ XX}^{1/2}f \bigr\rangle_{{\cal G}_{Y}} = \cov
\bigl[ g(Y), f(X) \bigr].
\end{equation}
%
\end{theorem}

The easiest way to understand equality (\ref{eqtarget1}) is through
the special case where
$f= f' - E (f') $, $g = g' - E (g')$ where $f' \in{\cal G}_{X}$, $g'
\in{{\cal G}_{Y}}$.
In this case,
\begin{eqnarray*}
\bigl\langle\tilde\Sigma_{ YY}^{1/2}g, \ryx\tilde
\Sigma_{ XX}^{1/2}f \bigr\rangle_{{\cal G}_{Y}} &=& \bigl\langle
\Sigma_{ YY}^{1/2}g', \ryx\Sigma_{ XX}^{1/2}h'
\bigr\rangle_{{\cal G}_{Y}} = \bigl\langle g', \coyx
f' \bigr\rangle_{{\cal G}_{Y}} \\
&=& \cov\bigl[f(X), g(Y)\bigr].
\end{eqnarray*}
The theorem also implies that, for all $f, g \in{L_2(P_{X})}$ and $s,
t \in
{L_2(P_{Y})}$,
\begin{eqnarray*}
\bigl\langle\tilde\Sigma_{ XX}^{1/2}g, \tilde
\Sigma_{ XX}^{1/2}f \bigr\rangle_{{\cal G}_{X}} &=& \langle g, f
\rangle_{{L_2(P_{X})}} = \cov\bigl[ g(X), f(X) \bigr],
\\
\bigl\langle\tilde\Sigma_{ YY}^{1/2}s, \tilde
\Sigma_{ YY}^{1/2}t \bigr\rangle_{{\cal G}_{Y}} &=& \langle s, t
\rangle_{{L_2(P_{Y})}} = \cov\bigl[ s(Y), t(Y) \bigr].
\end{eqnarray*}

%
%


\subsection{Generalized SIR}\label{sec5.2} 

The results of the last subsection allow us to characterize ${L_2 (P_{ X})}
\ominus{L_2 (P_{ Y})}$ in terms of
extended covariance operators, which is the key to developing its estimator.
Recall that classical SIR [\citet{Li91N2}] for linear SDR is based on the matrix
%
\begin{equation}
\label{eqclassicalsir} \bigl[\var(X)\bigr]^{-1}\var\bigl[ E (X | Y)
\bigr].
\end{equation}
Under the linear conditional mean assumption
requiring that $E(X|\beta^{\mathsf{T}}X)$ be linear in $X$ for any matrix
$\beta$ spanning
$\spc_{ Y|X}$, the re-scaled ``inverse'' conditional mean $[\var
(X)]^{-1}E(X|Y)$ is contained in this
space. To generalize this to the nonlinear setting, we first introduce
a conditional mean operator.

\begin{definition}\label{definitionexy}
We call the operator $\ecoyy^{-1/2}\ryx\ecoxx^{1/2}\dvtx  {L_2 (P_{
X})}\to {L_2 (P_{ Y})}$ the conditional expectation operator, and
denote it by $\exy$.
\end{definition}

The relation between the conditional expectation operator and
conditional expectations is
elucidated by the next proposition, which is followed by an important corollary.

\begin{proposition}\label{proposition2items} Under conditions \textup{(A)} and
\textup{(B)}, we have:
\begin{longlist}[(2)]
\item[(1)] for any $f \in{L_2 (P_{ X})}$,
$\exy f = E (f (X) | Y);$
\item[(2)] for any $g \in{L_2 (P_{ Y})}$, $\exy^*g = E ( g (Y) |X)$.\vadjust{\goodbreak}
\end{longlist}
\end{proposition}

\begin{pf}
For any $g \in{L_2 (P_{ Y})}$,
\begin{eqnarray*}
\langle\exy f, g \rangle_{L_2 (P_{ Y})}&=& \bigl\langle\ecoyy^{-1/2}\ryx
\ecoxx^{1/2}f, g \bigr\rangle_{L_2 (P_{ Y})} = \bigl\langle\ryx
\ecoxx^{1/2}f, \ecoyy^{1/2}g \bigr\rangle_\hyquo\\
&=& \cov
\bigl(f(X),g(Y)\bigr),
\end{eqnarray*}
where the last equality follows from (\ref{eqtarget1}). Hence
$\cov( f(X) - (\exy f)(Y),\break g(Y) ) = 0$.
By the definition\vspace*{1pt} of conditional expectation, $\exy f = E(f(X)|Y)$,
which proves 1. Assertion 2
follows from the fact that $\ecoyy^{-1/2}$ and $\ecoxx^{1/2}$ are
isomorphisms, and $\ryx^*= \rxy$.
\end{pf}

\begin{corollary} Under conditions \textup{(A)} and \textup{(B)}, for any $f, g \in
{L_2 (P_{ X})}$,
%
\begin{equation}
\label{eqcee} \bigl\langle g, \exy^*\exy f \bigr\rangle_{L_2 (P_{ X})}= \cov
\bigl[ E\bigl(g (X) |Y\bigr), E\bigl(f(X) | Y\bigr)\bigr].
\end{equation}
Moreover, $\exy^*\exy\in\B({L_2 (P_{ X})})$, and its norm is no
greater than 1.
\end{corollary}

\begin{pf}
We have
\begin{eqnarray*}
\bigl\langle g, \exy^*\exy f \bigr\rangle_{L_2 (P_{ X})}&=& \langle\exy g, \exy f
\rangle_{L_2 (P_{ X})} \\
&=& \bigl\langle E\bigl(g (X) |Y\bigr), E\bigl(f (X) |Y\bigr)
\bigr\rangle_{L_2 (P_{ X})},
\end{eqnarray*}
which is the right-hand side of (\ref{eqcee}). Moreover, since
$\ecoxx^{1/2}$ is isomorphic, we have
\[
\exy^*\exy= \bigl(\ecoyy^{-1/2}\ryx\ecoxx^{1/2}\bigr)^*\bigl(
\ecoyy^{-1/2} \ryx \ecoxx^{1/2}\bigr) = \ecoxx^{-1/2}\rxy
\ryx\tilde\Sigma_{ XX}^{1/2}.
\]
Hence
$
\| \exy^*\exy\| \le\| \ecoxx^{-1/2}\|  \| \rxy\|  \| \ryx\|
\| \ecoxx^{1/2}\|.
$
Because $\ecoxx^{1/2}$ and\break $\ecoxx^{-1/2}$ are isomorphisms, their norms
are both 1. By Baker [(\citeyear{Bak73}),
Theorem 1], $\| \ryx\| \le1$. Hence $\| \exy^*\exy\| \le1$.
\end{pf}

From this corollary we see that the quadratic form
\[
f \mapsto\bigl\langle f, \exy^*\exy f \bigr\rangle_{L_2 (P_{
X})},\qquad
{L_2 (P_{ X})} \times{L_2 (P_{ X})}
\to\real
\]
generalizes the matrix $\var[ E (X|Y)]$ of the linear case, which is
the essential ingredient of SIR for linear SDR. It is then not
surprising that the operator $\exy^*\exy$ is
closely connected to the central class for nonlinear SDR, as shown in
the following theorem. 

\begin{theorem}\label{theoremfisherforgsir} If conditions \textup{(A)} and
\textup{(B)} are satisfied and $\central$ is complete, then
\[
\overran\bigl(\exy^*\exy\bigr) = \central.
\]
\end{theorem}

\begin{pf}
By Lemma~\ref{lemmakeylemma}, $f \in\complete$ if and only if
$f \in{L_2(P_{X})}$ and $E(f |Y) = 0$. By
Proposition~\ref{proposition2items}, this happens
if and only if $f \in\ker\exy$. This shows\vadjust{\goodbreak} $\ker\exy= {\mathfrak
{C}_{ Y|X}^\perp}$.
However, because $\ker( \exy) = \ker( \exy^*\exy)$, we have
\[
\overran\bigl( \exy^*\exy\bigr) = \bigl[\ker\bigl( \exy^*\exy\bigr)
\bigr]^\perp= ( \ker\exy )^\perp= \bigl({\mathfrak{C}_{ Y|X}^\perp}
\bigr)^\perp= \complete.
\]
Since $\central$ is complete, we have $\complete= \central$, as
desired.
\end{pf}



Note that, unlike in classical SIR for linear SDR, here we do not have
to consider an analogue to the rescaling $[\var(X)]^{-1}$
in (\ref{eqclassicalsir}). This is because the ${L_2 (P_{ X})}$-inner
product absorbs the marginal variance
in the predictor vector. We refer to the sample
estimator based on $\overran(\exy^*\exy)$ (see Section
\ref{subsectionalgorithmofgsir}) as
\textit{generalized SIR} or GSIR.
The GSIR estimator is related to kernel canonical component analysis
(KCCA) introduced by \citet{BacJor03}; see also \citet{FukBacGre07}.
In Section~\ref{subsectionalgorithmofgsir} we will explore
similarities and differences between
these two methods.

\subsection{Kernel SIR}\label{sec5.3}

We now turn to another nonlinear SDR estimator, which was proposed by
\citet{Wu08} and further
studied by \citet{YehHuaLee09}, called \textit{kernel sliced inverse
regression}
(KSIR).
In our setting,
the population-level description of this estimator is as follows. Let
$\hx$ be a Hilbert space
satisfying (A) and (B) (in this case an RKHS, but this assumption is
unnecessary). Let
$T\dvtx  \hx\to{L_2(P_{X})}$ be the centering transformation $T(f) = f - E(f)$.
Let $J_1,\ldots, J_h$ be a partition of
$\oy$, and let $\mu_1,\ldots, \mu_h \in\overran T$ be the Riesz
representations of the linear
functionals
\[
T_j\dvtx \overran T \to\real,\qquad g \mapsto E\bigl( g (X) | Y \in
J_i \bigr),\qquad i=1,\ldots, h.
\]
In our language, KSIR uses (the sample version of) the subspace $\sp
(\coxx^{-1}\mu_1,\break\ldots,
\coxx^{-1}\mu_h )$ to estimate $\central$. The next theorem shows that
any such representation
must be a member of
$\complete$, and thus of $\central$ (since $\complete\subseteq
\central
$)---which implies that
KSIR is unbiased.

\begin{theorem} If \textup{(A)} and \textup{(B)} hold, then $\mu_j \in\complete$.
\end{theorem}

\begin{pf}
By condition (A), $\overran T = {L_2(P_{X})}$.
If $f \in{L_2 (P_{ X})}\ominus{L_2 (P_{ Y})}\subseteq\overran T$, then, by
Lemma~\ref{lemmakeylemma},
$E(f | Y) = 0 $. Hence
$
\langle f, \mu_i \rangle_{L_2 (P_{ X})}= E [f(X) |Y\in J_i ] = 0
$.
\end{pf}

\citet{YehHuaLee09} give another unbiasedness proof for KSIR, but
they assume that
the spanning functions of $\hx$, say $f_1,\ldots, f_m$, satisfy the
linear conditional mean
assumption. That is, for any $f \in\hx$,
$E(f|f_1,\ldots, f_m)$ has the form $c_0 + c_1 f_1 + \cdots+ c_m f_m$
for some $c_0,\ldots, c_m \in\real$.
This condition is an analogue of the linear conditional mean assumption
for linear SDR; see,
for example, \citet{Li91N2} and \citet{CooLi02}. Interestingly, our
result no longer relies on
this assumption. The reason Yeh, Huang and Lee need the assumption in
the first place is that they
define the central class [Definition 1 of \citet{YehHuaLee09}] as
the linear subspace spanned by $h_1,\ldots, h_d$ in $\sp(f_1,\ldots,
f_m)$ such that
%
\begin{equation}
\label{eqhh} Y \indep X | h_1(X),\ldots, h_d(X),
\end{equation}
whereas we define the central class as the class of all measurable
functions of $h_1,\ldots, h_d$.
Indeed, in the nonlinear setting there is no reason to restrict to this
linear span formulation,
since the conditional independence (\ref{eqhh}) only relies on the
$\sigma$-field generated by
$h_1,\ldots, h_d$.

\section{Beyond completeness: Generalized SAVE}\label{sec6}\label{sectionbeyond}

We now turn to the more general problem of estimating the central class
when it is not complete,
in which case the regression class may be a proper subset of the
central class.
We will generalize SAVE [\citet{Li91N1}] to the nonlinear
case and show that it can recover
functions beyond the regression class.

The setting here is different from that for GSIR in two respects.
First, since we now deal with
the location-invariant quantity $f(X) - E[f(X)|Y]$, we no longer need
to define the conditional mean
operator through the centered $L_2$-spaces ${L_2 (P_{ Y})}$ and ${L_2
(P_{ X})}$. Second, we
now define relevant operators
through $L_2$-spaces instead of RKHSs, which is more convenient in this context.
Let ${{L_2' (\px)}}$ and ${{L_2' (\py)}}$ denote the noncentered
$L_2$-spaces. Define the
noncentered conditional mean operator $E_{ X|Y}'\dvtx  {{L_2' (\px)}}\to
{{L_2' (\py)}}$ through
%
\begin{equation}
\label{eqnewe}\quad \bigl\langle g, E_{ X|Y}'f \bigr
\rangle_{L_2'(P_{ X})} = E\bigl(g(Y) f(X)\bigr),\qquad f\in{{L_2'
(\px)}}, g \in{{L_2' (\py)}}.
\end{equation}
By the same argument of Proposition~\ref{proposition2items}, $E_{ X|Y}'f
= E(f(X)|Y)$.
To generalize SAVE, we introduce a new type of conditional variance operator.

\begin{definition}\label{definitionhetcon}
For each $y \in\oy$, the bilinear form
\[
{L_2(P_{X})}\times{L_2(P_{X})}\to
\real,\qquad (f, g ) \mapsto\bigl(E_{
X|Y}'(f g ) -
E_{ X|Y}'f E_{ X|Y}'g\bigr) (y)
\]
uniquely defines an operator $\vxy(y) \in\B({L_2(P_{X})})$ via the Riesz
representation. We call the random
operator
\[
\vxy\dvtx \oy\to\B\bigl({L_2(P_{X})}\bigr),\qquad y \mapsto
\vxy(y)
\]
the heteroscedastic conditional variance operator given $Y$.
\end{definition}
The operator $\vxy$ is different from the conditional variance operator
$\cxy$ introduced by Fukumizu,
Bach and Jordan (\citeyear{FukBacJor03}, \citeyear{FukBacJor09}). In a sense, $\cxy$ is a generalization of
$E[\var(X|Y)]$ rather than
$\var(X|Y)$, because $\langle f, \cxy f \rangle_\hx= E[\var(f(X)|Y)]$.
Note that $E[\var(f(X)|Y)]$
becomes\break $\var(f(X)|Y)$ only when the latter is nonrandom. So $\cxy$
might be called a \textit{homoscedastic} conditional variance operator.\vadjust{\goodbreak} In
contrast, $\langle f, \vxy f \rangle_{L_2(P_{X})}$ gives directly the
conditional variance $\var[f(X)|Y]$, hence the term heteroscedastic
conditional variance operator.
Here, we should also stress that $E_{ X|Y}'$ is defined between noncentered
${{L_2' (\px)}}$ and ${{L_2' (\py)}}$,
whereas $\vxy(y)$ is defined between centered ${L_2 (P_{ X})}$ and
${L_2 (P_{ X})}$.

We now define the expectation of a generic random operator $A\dvtx  \oy\to
\B({L_2(P_{X})})$. For each
$f \in{L_2(P_{X})}$ and $x \in\ox$, the mapping $y \mapsto(A(y)f)(x)$
defines a random variable. Its
expectation defines a function $x \mapsto\break\int_{\oy} (A(y)f)(x) \*\py
(dy)$, which is a member of ${L_2(P_{X})}$. Denoting this member as
$\tilde
f$, we define the nonrandom operator ${L_2(P_{X})}\to{L_2(P_{X})},
f\mapsto
\tilde f$ as the expectation $E (A)$. We now consider the operator
%
\begin{equation}
\label{eqsaveop} S = E ( V - \vxy)^2\dvtx {L_2(P_{X})}
\to{L_2(P_{X})},
\end{equation}
where $V\dvtx  {L_2(P_{X})}\to{L_2(P_{X})}$ is the (unconditional)
covariance operator
defined by
\[
\langle f, V g \rangle_{L_2 (P_{ X})}= \cov\bigl(f(X), g(X)\bigr).
\]
This operator is similar to $\ecoxx$ in Section~\ref{sectiongsir}
except that it is not defined through
RKHS. The operator $S$ is an extension of the SAVE matrix [\citet{Li91N1}]
%
\begin{equation}
\label{eqsavemat} \Sigma^{-1}E \bigl[ \var(X) - \var(X|Y)
\bigr]^2 \Sigma^{-1}.
\end{equation}
%
%
Let $\beta$ be a basis matrix of the central subspace $\lcentral$ of
linear SDR. Cook and Weisberg show
that if $E(X|\beta^{\mathsf{T}}X)$ is linear in $\beta^{\mathsf
{T}}X$ and
$\var(X | \beta^{\mathsf{T}}X)$ is nonrandom, then the column space of
(\ref
{eqsavemat}) is contained
in $\lcentral$. The next theorem generalizes this result, but without
requiring an analogue of the linear
conditional mean assumption.


\begin{theorem}\label{theoremsave1} Suppose that conditions \textup{(A)} and
\textup{(B)} are satisfied, and
$\var[f(X)|\cfield]$ is nonrandom for any $f \in{\mathfrak S}_{
Y|X}^{\perp}$. Then $\overran S \subseteq\central$.
\end{theorem}

\begin{pf}
Let $f \perp\central$.
We claim that for any
$y \in\oy$,
%
\begin{equation}
\label{eqh,f=0} \bigl\langle f, \bigl[V - \vxy(y)\bigr] f \bigr
\rangle_{L_2(P_{X})}= 0.
\end{equation}
%
Because $Y \indep X | \cfield$, we have
\[
\var\bigl( f(X) | Y\bigr) = \var\bigl( E \bigl(f(X) | \cfield\bigr) | Y\bigr) +
E \bigl(\var\bigl( f(X) | \cfield\bigr) | Y \bigr).
\]
%
Because, by Lemma~\ref{lemmakeylemma}, $E(f(X) | \cfield)$ is
constant, the first term is 0.
Because $\var( f(X) | \cfield)$ is nonrandom, the second term is
$\var
(f(X) | \cfield)$. Hence
\[
\var\bigl( f(X) | Y\bigr) = \var\bigl( f(X) | \cfield\bigr).
\]
Similarly,
\[
\var\bigl( f(X) \bigr) = \var\bigl( E \bigl( f(X) | \cfield\bigr) \bigr) + E
\bigl( \var\bigl(f(X) | \cfield\bigr) \bigr) = \var\bigl(f(X) | \cfield\bigr).\vadjust{\goodbreak}
\]
Therefore $\var( f(X) | Y) = \var( f(X))$, which implies (\ref
{eqh,f=0}). Since $V - \vxy(y)$ is self-adjoint, (\ref{eqh,f=0})
implies $f \in\ker\vxy(y)$.
Hence
\[
\bigl\langle f, \bigl[ V - \vxy(y)\bigr]^2 f \bigr
\rangle_{L_2(P_{X})}= 0.
\]
Now integrate both sides of this equation to obtain
\begin{eqnarray*}
&&
\int_{\oy} \bigl\langle f, \bigl(V - \vxy(y)
\bigr)^{2}f \bigr\rangle_{L_2(P_{X})}\py(dy) \\
&&\qquad= \biggl\langle f,
\int_{\oy}\bigl(V - \vxy(y)\bigr)^{2}f\py(dy)
\biggr\rangle_{L_2(P_{X})}
\\
&&\qquad= \bigl\langle f, \bigl(E (V-\vxy)^{2}\bigr) f \bigr
\rangle_{L_2(P_{X})}= 0.
\end{eqnarray*}
Hence $f \in\ker E (V - \vxy)^2$, as desired.
\end{pf}

Similar to the case of GSIR, we do not need to employ the rescaling
by $\Sigma^{-1}$ in (\ref{eqsavemat}) when generalizing
SAVE, because
the ${L_2(P_{X})}$-inner product absorbs any marginal variance. We call the
estimator derived from $\overran S$
(see Section~\ref{subsectionalgorithmofgsave}) \textit{generalized
SAVE} or GSAVE.
The next theorem shows that GSAVE can recover functions outside
$\complete$.

\begin{theorem}\label{theoremsave2} If conditions \textup{(A)} and \textup{(B)} are
satisfied, then
$
\complete\subseteq\overran S.
$
\end{theorem}

\begin{pf}
Since $S$ is self-adjoint, it suffices to show that $\ker S
\subseteq{\mathfrak{C}_{ Y|X}^\perp}$.
For any $f \in\ker S$,
\[
\int_{\oy} \bigl\langle f, (V - \vxy)^{2}(y) f
\bigr\rangle\py(dy) = 0.
\]
Hence $\langle f, ( V - \vxy(y))^{2}f \rangle_{L_2(P_{X})}= 0$ a.s.
$\py
$, which implies
$(V - \vxy(y)) f = 0$ a.s. $\py$. Then
\[
\int_{\oy} \bigl\langle f, \bigl(V - \vxy(y)\bigr) f \bigr
\rangle_{L_2(P_{X})}\py(dy) = 0.
\]
By Definition~\ref{definitionhetcon}, the left-hand side is $\var
[f(X)] - E [\var(f(X)|Y)]
=\break \var[E(f(X)| Y)]$. Hence $\var[ E(f(X) | Y)] = 0$, which implies
$E[f(X)|Y]= E[f(X)] = 0$.
By Lemma~\ref{lemmakeylemma}, we have $f \in{L_2 (P_{ X})}\ominus
{L_2 (P_{ Y})}=
{\mathfrak{C}_{ Y|X}^\perp}$, as desired.
\end{pf}


Combining Theorems~\ref{theoremsave1} and~\ref{theoremsave2} we see that
%
\begin{equation}
\label{eqsave3} \complete\subseteq\overran S \subseteq\central,
\end{equation}
which is analogous to the relation
(\ref{eqcookcritchley}) in the classical setting.
Thus we can expect GSAVE to discover functions outside the class
$\complete$, just as we can expect SAVE
to discover vectors outside the space spanned by SIR.

\section{Algorithms}\label{sec7}

We now develop algorithms for the sample versions of GSIR and GSAVE,
together with a cross-validation scheme
to select parameters in the GSIR and GSAVE algorithms.
These sample versions involve representing the operators
in Theorems~\ref{theoremfisherforgsir} and~\ref{theoremsave1} as matrices.
To formulate the algorithms we need to introduce coordinate
representations of functions and operators,
which we adopt with modifications
from Horn and Johnson [(\citeyear{HorJoh85}), page 31]; see also \citet{LiChuZha}.

Throughout this section, $A^{\dagger}$ represents the Moore--Penrose
inverse of a matrix $A$, $A^{{\dagger} \alpha}$ represents
$(A^{\dagger} )^\alpha$, $I_n$ denotes the $n\times n$ identity matrix,
$1_n$ denotes the vector in $\real^n$ whose entries are all 1 and
$Q=I_n - 1_n 1_n^{\mathsf{T}}/n$. Let $\kax\dvtx  \ox\times\ox\to\real$
be a positive definite function. Also, let $\kx$ be the $n \times n$
the Gram matrix $\{\kax(X_i,X_j)\dvtx  i,j = 1,\ldots, n \}$, $\gx$ its
centered versions $Q \kx Q$ and $\llx$ the Gram matrices with
intercept; that is, $\llx= ( 1_n, \kx)^{\mathsf{T}}$. Finally, define
$\kay, \kym, \gy, \lly$ in the same manner for $Y$.

\subsection{Coordinate representation}\label{sec7.1}

Let $\H$ be a finite-dimensional Hilbert space with spanning system
$\B
= \{b_1,\ldots, b_n \}$.
For an $f \in\H$, let $[f]_{\B} \in\real^n$ denote the coordinates of
$f$ relative to $\B$; that is,
$f = \sum_{i=1}^n ([f]_{\B})_i b_i$. Let $b\dvtx \ox\to\real^n$ denote the
$\real^n$-valued function
$(b_1,\ldots, b_n )^{\mathsf{T}}$. Then we can write $f = [f]_\B^{\mathsf{T}}b$.
Let $A\dvtx  \H\to\H'$, where $\H'$ is another finite-dimensional Hilbert
spaces with spanning system
$\C= \{c_1,\ldots, c_m\}$ and let $c = (c_1,\ldots, c_m)^{\mathsf{T}}$.
Then, for $f \in\H$,
\[
Af = A \bigl( b^{\mathsf{T}}[f]_{\B} \bigr) = (A b_1,
\ldots, A b_n) [f]_{\B
} = \bigl( c^{\mathsf{T}}[Ab_1]_{\C},
\ldots, c^{\mathsf{T}}[A b_n]_{\C
} \bigr)
[f]_{\B}.
\]
Thus, if we let ${_\C}[ A ]_{\B} = ( [ A b_1 ]_{\C},\ldots, [A
b_n]_{\C
} )$, then
$A f = c^{\mathsf{T}}({_\C}[A]_{\B}) [f]_{\B}.$ In other words,
\[
[A f]_{\C} = \bigl( {_\C}[A]_{\B}\bigr)
[f]_{\B}.
\]
Furthermore, if $A_1\dvtx  \H' \to\H''$ is another linear operator, where
$\H''$ is a third finite-dimensional
Hilbert space with spanning system $\D$, then, by a similar argument,
\[
{_\D}[A_1A]_{\B} = \bigl({_\D}[A_1]_{\C}
\bigr) \bigl({_\C}[A]_{\B}\bigr).
\]
Since the spanning systems in the domain and range of an operator are
self-evident in the following discussion,
we will write ${_{ {\cal C}}}[A]_{ {\cal B}}$ and $[f]_{ {\cal B}}$
simply as
$[A]$ and $[f]$. 

Suppose $A \in\B(\H)$ is self-adjoint.
It can be shown that, for any $\alpha> 0$, \mbox{$[A^\alpha] = [A]^\alpha$}.
Depending on the choice of the spanning system of $\H$, it is possible
that $A$ is invertible
and yet $[A]$ is singular, but it is generally true that $A^{-\alpha} =
[A]^{\dagger\alpha}$. Throughout this section the square brackets
\mbox{$[\cdot]$} will be used exclusively for denoting
coordinate representations.

%

\subsection{Algorithm for GSIR}\label{sec7.2}\label{subsectionalgorithmofgsir}


At the sample level, $P_{ X}$ is replaced by
the empirical measure $P_{ n,X}$; $\hx$ is the RKHS spanned by
$\bx= \{\kax(\cdot, X_1),\ldots,\break \kax(\cdot,X_n)\}$ with inner
product $\langle f, g \rangle_{\hx} = [f]^{\mathsf{T}}\kx[g]$,
where $[\cdot]$ is coordinate
with respect to~$\bx$.
The space $L_2 (P_{ n,X})$ is spanned by
$ \kax(\cdot, X_i) - E_n \kax(X, X_i)$, $i=1,\ldots, n$,
with inner product $\langle f, g \rangle_{ L_2 (P_{n,X})} = \cov_n
[f(X),\break g(X)] = n^{-1}[f] \kx Q \kx[g]$.
The operator
$\mxx$ is defined through the relation $\langle f, \mxx g\rangle_{
\hx} = \cov_n (f(X), g(X))$; that is,
\[
[f]^{\mathsf{T}}\kx[\mxx] [g] = n^{-1}[f]^{\mathsf{T}}\kx Q \kx[g].
\]
Since $[f]$ and $[g]$ are arbitrary members of $\real^n$, the above
implies $[\mxx] = n^{-1}Q \kx$. Then any $f \in\ran\mxx\equiv
{\cal
G}_{ X}$ can be written as
$ \mxx g$ for some $g \in\hx$, which implies $[f] = Q \kx[g] = Q[f]$.
Consequently, for any $f, g \in\gfam$,
$\langle f, g \rangle_{\hx} = [f]^{\mathsf{T}}\gx[g]$.


Let us now find the matrix representations of $\coxx$, $\coyy$ and
$\coyx$.
In the following, $\hbasis$ represents the function $x \mapsto(\kax
(x,X_1),\ldots, \kax(x, X_n))^{\mathsf{T}}$.
For any $f \in\gclass$, we have
\[
\coxx f = \mxx f = h_{ X}^{\mathsf{T}}[\mxx] [f] =
n^{-1}h_{
X}^{\mathsf{T}}Q \kx Q [f] = n^{-1}h_{ X}^{\mathsf{T}}
\gx[f].
\]
Hence $[\coxx f] = [\coxx][f] = n^{-1}\gx[f]$.
Since this is true for all $[f] \in\sp(Q)$, we have $[\coxx] = n^{-1}
\gx$.
By the same argument we can show that
%
\begin{eqnarray}
\label{eqrepresentation}
[\ecoxx] &=& n^{-1}\gx,\qquad [\coyy] = [\ecoyy] =
n^{-1}\gy,\nonumber\\[-8pt]\\[-8pt]
[\coyx] &=& [\ecoyx] = n^{-1}\gx,\qquad
[\exy] = \gy^{\dagger}\gx\gx^{{\dagger1/2}} \gx^{{1/2}}.
\nonumber
\end{eqnarray}

Theorem~\ref{theoremfisherforgsir} suggests that we use
$\overran( \exy^*\exy)$ to estimate $\central$. Since $\exy^*
\exy$
is an operator on ${L_2 (\pnx)}$ to ${L_2 (\pnx)}$,
the vectors in\break $\overran( \exy^*\exy)$ can be found by
\[
\mbox{maximizing}\quad \bigl\langle f, \exy^*\exy f \bigr\rangle_{
{L_2 (\pnx)}}
= \| \exy f \|_{{L_2(P_{n, { Y}})}}^{2}
\]
subject to
\[
\langle f, f\rangle_{{L_2 (\pnx)}} = 1.
\]
The coordinate representation of this problem is
\[
\mbox{maximizing}\quad [f]^{\mathsf{T}}[ \exy]^{\mathsf{T}}\gy^{2}
[ \exy] [f]\quad\mbox{subject to}\quad [f]^{\mathsf{T}}\gx^{2} [f] =
1.
\]
The optimal solution is $[f] = \gx^{{\dagger}} \phi$, where
$\phi$
are the leading eigenvectors of the matrix
%
\begin{eqnarray}
\label{eqgsirmatrix}
&&
\gx^{{\dagger}} [ \exy]^{\mathsf{T}}\gy^{2} [
\exy] \gx^{{\dagger}}\nonumber\\[-8pt]\\[-8pt]
&&\qquad= \gx^{{\dagger}} \gx^{1/2}\gx^{{\dagger1/2}}
\gx\gy^{{\dagger}} \gy^2 \gy^{{\dagger}} \gx\gx^{{\dagger1/2}}
\gx^{1/2} \gx^{{\dagger}}.\nonumber
\end{eqnarray}
To enhance accuracy we replace the Moore--Penrose inverses $\gx^{
\dagger}$ and $\gy^{\dagger}$ by the ridge-regression-type
regularized inverses $(\gx+ \epsx I_n)^{-1}$ and $(\gy+ \epsy
I_n)^{-1}$.
We summarize the algorithm as follows:
\begin{longlist}[(3)]
\item[(1)] Select the parameters $\gamx$, $\gamy$, $\epsx$, $\epsy$
using the algorithm in Section~\ref{subsectioncrossvalidation}.
\item[(2)] Compute the matrix
\[
(\gx+ \epsx I_n)^{-3/2} \gx^{3/2} ( \gy+ \epsy
I_n)^{-1} \gy^{2} ( \gy+ \epsy
I_n)^{-1} \gx^{3/2} (\gx+ \epsx
I_n)^{-3/2}
\]
and its first $d$ eigenvectors $\phi_1,\ldots, \phi_d$ of this matrix.
\item[(3)] Form the sufficient predictors at $x$
$
\phi_i^{\mathsf{T}}( \gx+ \epsx I_n)^{-1}h_{ X} (x)$, $i=1,\ldots, d$.
\end{longlist}

GSIR estimation is similar to the kernel canonical correlation analysis
(KCCA) developed by \citet{Aka01}, \citet{BacJor03} and
\citet{FukBacGre07}. In our notation, KCCA maximizes
\[
\langle g, \coyx f \rangle_{L_2 (P_{ Y})}= [g ]^{\mathsf{T}}\gy\gx [f ]
\]
subject\vspace*{1pt} to $\langle g, \coyy g \rangle_{L_2 (P_{ Y})}= [g ]^{\mathsf
{T}}\gy^{2}
[f ] = 1$
and $\langle f, \coxx f \rangle_{L_2 (P_{ X})}=[ g ]^{\mathsf{T}}\*\gy^{2} [g ]=
1$. The optimal solution
for $[f]$ is $[f ] = (\gx+ \epsilon I_n)^{-1}\phi$, where $\phi$ is
one of the first $d$ eigenvectors of
%
\[
(\gx+ \epsilon I_n)^{-1}\gx\gy(\gy+ \epsilon
I_n)^{-2} \gy \gx(\gx+ \epsilon I_n)^{-1}.
\]
%
We will compare GSIR and KCCA in Section~\ref{sec8}.

\subsection{Algorithm for GSAVE}\label{sec7.3}\label{subsectionalgorithmofgsave}

We first derive the sample-level representation of the operator $\vxy(y)$.
The sample version of the noncentered $L_2$-classes ${L_2' (P_{n,{
X}})}$ and
${L_2' (P_{n,{ Y}})}$ are spanned by
%
\begin{equation}
\label{eqspanningsystems} {\cal C}_{ X} = \bigl\{1, \kax(\cdot,
X_1),\ldots, \kax(\cdot, X_n)\bigr\},\qquad {\cal
C}_{ Y} = \bigl\{1, \kay(\cdot, Y_1),\ldots, \kay(\cdot,
Y_n)\bigr\},\hspace*{-35pt}
\end{equation}
respectively.
Let $[\cdot]$ represent the coordinates relative to these spanning systems.
Then, for any $f \in{L_2' (P_{n,{ X}})}$, $(f(X_1),\ldots, f
(X_n))^{\mathsf{T}}
=L_{ X}^{\mathsf{T}}[f]$. 
The operator $E_{ X|Y}'$ is defined through the relation
$
\langle g, E_{ X|Y}'f \rangle_{L_2' (P_{n,{ Y}})}= E_n (g(Y) f(X))
$, which yields the
representation
%
\begin{equation}
\label{eqexypco} \bigl[ E_{ X|Y}'\bigr] = \bigl(\lly
\lly^{\mathsf{T}}\bigr)^{\dagger}\bigl(\lly\llx^{\mathsf{T}}\bigr).
\end{equation}
Let $\llly$ denote the function $y \mapsto(1, \kay(y, Y_1),\ldots,
\kay(y, Y_n))^{\mathsf{T}}$, and let $\lllx$ denote the
same function of $x$. For any $f, g \in{L_2' (P_{n,{ X}})}$,
%
\begin{eqnarray}
\label{eqconvar} &&\bigl\{E_{ X|Y}'(fg)-
\bigl(E_{ X|Y}'f\bigr) \bigl(E_{ X|Y}'g
\bigr)\bigr\}(y)
\nonumber\\[-8pt]\\[-8pt]
&&\qquad= \llly^{\mathsf{T}}(y) \bigl[ E_{ X|Y}'\bigr]
[fg] - [f]^{\mathsf{T}}\bigl[E_{ X|Y}'
\bigr]^{\mathsf{T}}\llly(y) \llly^{\mathsf{T}}(y) \bigl[E_{ X|Y}'
\bigr] [g].
\nonumber
\end{eqnarray}
For any $X_i$, $f(X_i) g(X_i)$ can be expressed as the $i$th entry of
the vector
$ \llx^{\mathsf{T}}[f] \odot\llx^{\mathsf{T}}[g]$,
which is the same as $\llx^{\mathsf{T}}(\llx\llx^{\mathsf
{T}})^{\dagger}\llx(\llx^{\mathsf{T}}[f] \odot\llx^{\mathsf{T}}[g])$, where $\odot$ is the
Hadamard product.
Thus we have the coordinate representation
%
\begin{equation}
\label{eqfgco} [fg] = \bigl(\llx\llx^{\mathsf{T}}\bigr)^{\dagger}\llx\bigl(
\llx^{\mathsf{T}}[f] \odot\llx^{\mathsf{T}}[g]\bigr).
\end{equation}
Substituting (\ref{eqexypco}) and (\ref{eqfgco}) into (\ref
{eqconvar}) we see that, for any $f, g \in L_2'(P_{ n,X})$,
%
\begin{eqnarray}
\label{eqvxynoncentered}
\bigl\langle f, \vxy(y) g \bigr\rangle_{L_2'(P_{ n,X})} &=&
[f]^{\mathsf{T}}\llx\bigl( \diag\cy(y) - \cy(y) \cy^{\mathsf{T}}(y) \bigr)
\llx^{\mathsf{T}}[g] \nonumber\\[-8pt]\\[-8pt]
&\equiv&[f]^{\mathsf{T}}\llx\Lambda(y) \llx^{\mathsf{T}}[g],\nonumber
\end{eqnarray}
where $\cy(y) = \lly^{\mathsf{T}}(\lly\lly^{\mathsf{T}})^{\dagger}
\llly(y)$.

Let $S_n\dvtx  {L_2 (P_{ n, X})}\to{L_2 (P_{ n, X})}$ be the operator $E_n
(V - \vxy(Y))^2$.
By Theorem~\ref{theoremsave1}, GSAVE is the class of functions $\overran
(S)$. At the sample level, this
corresponds to
%
\begin{equation}
\label{eqgsaveobj} \mbox{maximizing}\quad \langle f, S_n f
\rangle_{L_2 (P_{ n, X})} \quad\mbox{subject to}\quad \langle f, f
\rangle_{L_2 (P_{ n, X})}=1.
\end{equation}
By (\ref{eqvxynoncentered}), for each $y \in\oy$, and $f, g \in
{L_2 (\pnx)}$, we have
\[
\bigl\langle g, \vxy(y) f \bigr\rangle_{L_2 (P_{n,{ X}})}= [f]^{\mathsf
{T}}\llx Q
\Lambda(y) Q \llx^{\mathsf{T}}[g].
\]
From this we deduce that $[\vxy(y)] = (\llx Q \llx^{\mathsf{T}}/
n)^{\dagger}
\llx Q \Lambda(y) Q \llx^{\mathsf{T}}$. By a similar derivation
we find $[V] = (\llx Q \llx^{\mathsf{T}}/n)^{\dagger}(\llx Q \llx^{\mathsf{T}}/n)$. Hence
\[
\bigl[V - \vxy(y)\bigr] = \bigl(\llx Q \llx^{\mathsf{T}}/n\bigr)^{\dagger}
\llx Q \bigl(Q/n - \Lambda(y)\bigr) Q \llx^{\mathsf{T}}.
\]
It follows that
\[
\langle f, S_n f \rangle_{L_2 (P_{n,{ X}})} = E_n \bigl\{
[f]^{\mathsf{T}}\llx Q \bigl( Q / n - \Lambda(Y)\bigr) Q \bigl( Q / n -
\Lambda(Y)\bigr) Q \llx^{\mathsf{T}}[f] \bigr\}.
\]
To find $\ran(S_n)$ we maximize the above
subject to $[f]^{\mathsf{T}}(\llx Q \llx^{\mathsf{T}}/ n ) [f] = 1$.
Again we use the regularized inverses instead of the Moore--Penrose
inverses to enhance performance. 
The algorithm is summarized as follows:
\begin{longlist}[(4)]
\item[(1)] Determine $\gamx, \gamy, \epsx, \epsy$ using the algorithm
is Section~\ref{subsectioncrossvalidation}.
\item[(2)]
Compute
$
C = \elly^{\mathsf{T}}(\elly\elly^{\mathsf{T}}+ \epsy I_{n+1} )^{-1/2}
\elly$.
Let $C_i$ be the $i$th
column of $C$. Compute $\Lambda_i = \diag(C_i) - C_i C_i^{\mathsf
{T}}$ and
then compute $\Gamma_i = Q/n - \Lambda_i$
for $i=1,\ldots, n$.
\item[(3)] Compute
\[
n^{-1}\sum_{i=1}^n \bigl(\llx Q
\llx^{\mathsf{T}}+ \epsx I_{n+1} \bigr)^{-1/2} \llx Q
\Gamma_i Q \Gamma_i Q \llx^{\mathsf{T}}\bigl(\llx Q
\llx^{\mathsf{T}}+ \epsx I_{n+1} \bigr)^{-1/2}
\]
and the first $d$ eigenvectors of this matrix, say $\phi_1,\ldots,
\phi_d$.
\item[(4)] The sufficient predictors' values at any given $x \in\ox$
are the set of $d$ numbers
\[
\lllx^{\mathsf{T}}(x) \bigl(\llx Q \llx^{\mathsf{T}}+ \epsx I_{n+1}
\bigr)^{-1/2} Q \phi_i,\qquad i =1,\ldots, d.
\]
\end{longlist}

Here we should mention that, similar to SAVE for linear SDR, GSAVE
works best for extracting predictors
affecting the conditional variance of the response, but often not so
well for extracting predictors
affecting the conditional mean. However, we expect that other
second-order methods for linear SDR, such
as directional regression [\citet{LiWan07}] and the minimum
discrepancy approach [\citet{CooNi05}],
will be amenable to similar generalizations to nonlinear SDR. These
will be left for future research.

\subsection{Cross-validation algorithm}\label{sec7.4}\label{subsectioncrossvalidation}

We now develop a cross-validation scheme to determine the parameters
$\gamx$, $\gamy$, $\epsx$, $\epsy$, which are used in the algorithms
for both the GSIR and the GSAVE.
We will only describe the algorithm for determining $(\gamx, \epsx)$;
that of
$(\gamy, \epsy)$ is completely analogous.

In the following, for a matrix $A$, $A_{-i,-j}$ represents the
submatrix of $A$ with its
$i$th row and $j$th column removed, and $A_{-i,j}$ represents the
$j$th column of $A$ with the $i$th entry removed.
Let ${\cal C}_{ Y}^{-i}={\sub{{\cal C}} Y} \setminus\{ \sub
\kappa
Y ( \cdot, Y_i) \}$, and
define ${\cal C}_{ X}^{-i}$ similarly.
Our cross-validation strategy is to predict $f(Y_i)$ for each $f \in
{\cal C}_{ Y}^{-i}$, using the conditional mean operator developed
from $({\cal C}_{ X}^{-i}, {\cal C}_{ Y}^{-i})$.
The regularized matrix representation of $\sub{E'} {Y|X}$ based on
$({\cal C}_{ X}^{-i}, {\cal C}_{ Y}^{-i})$
is
\begin{eqnarray*}
&&
\bigl[ \sub{E'} {Y|X} \bigr]_{-(i+1),-(i+1)} \\
&&\qquad= \bigl\{(\sub L X
)_{-(i+1),-i} (\sub L X )_{-(i+1),-i}^{\mathsf{T}}+ \sub\epsilon X
I_n \bigr\}^{-1}(\sub L X)_{-(i+1),-i} (\sub L
Y)_{-(i+1),-i}^{\mathsf{T}}.
\end{eqnarray*}
The $k$th member $f_k$ of ${\cal C}_{ Y}^{-i}$ is the function $
e_k^{\mathsf{T}}{(\sub
\ell Y)}_{-(i+1)} (\cdot) $
where $e_k$ is the vector in $\real^n$
whose $k$th entry is 1 and the remaining entries are 0. Therefore, the
estimate of $E(f_k(Y)|X=x)$ based on
on ${\cal C}_{ X}^{-i}$ is
\[
(\sub\ell X)_{-(i+1)}^{\mathsf{T}}(x) \bigl[ \sub{E'}
{Y|X} f_k \bigr]_{-(i+1)} = e_k^{\mathsf{T}}
\bigl[ \sub{E'} {Y|X}\bigr]_{-(i+1),-(i+1)}^{\mathsf{T}}(\sub \ell
X)_{-(i+1)} (x),
\]
and the prediction of $(f_1(Y_i),\ldots, f_n (Y_i))^{\mathsf{T}}$ is
$[ \sub{E'}{Y|X}]_{-(i+1),-(i+1)}^{\mathsf{T}}(\sub\ell X)_{-(i+1)} (X_i)$.
However, because $(\sub\ell X)_{-(i+1)} (X_i)$
is the vector $(\sub L X)_{-(i+1), i}$, and $(f_1(Y_i),\ldots,\break f_n
(Y_i))^{\mathsf{T}}$
is the vector $(\sub L Y)_{-(i+1), i}$, the difference between
$(f_1(Y_i),\ldots, f_n (Y_i))^{\mathsf{T}}$
and its prediction is
\[
(\sub L Y)_{-(i+1), i} - \bigl[ \sub{E'} {Y|X}
\bigr]_{-(i+1),-(i+1)}^{\mathsf
{T}}(\sub L X)_{-(i+1), i}.
\]
To stress that this difference depends on $\gamx, \epsx, \gamy$, we
denote it by
$\Delta_i (\epsx,\break \gamx, \gamy)$.
Our cross-validation criterion is defined as $
\cv(\gamma_{ X}, \epsilon_{ X}, \gamy) = \break\sum_{i=1}^n \|
\Delta_i ( \gamma_{ X}, \epsilon_{ X}, \gamma_{ Y} )
\|^2
$. Since the role of $\gamy$
is only to determine the set of functions to be predicted, we exclude
it from the optimization process (for the determination of $\epsx,
\gamx
$). Moreover,
as argued in \citet{FukBacJor09}, the parameters $\epsx$
and $\gamx$ have similar smoothing
effects and only one of them needs to be optimized. For these reasons
we fix $\gamy$ and $\epsx$ at
%
\begin{equation}
\label{eqinitialvalue} 1/\gamma_{ Y 0} = \pmatrix{n
\cr
2}^{-1}
\sum_{i < j } | Y_i - Y_j
|^2,\qquad \epsilon_{X 0} = 0.01
\end{equation}
and minimize $\cv(\gamx, \epsilon_{X 0}, \gamma_{Y 0})$
over a
grid for $\gamx$.
The grid consists of 20 subintervals in
$[\gamma_{X 0}/3, 3 \gamma_{X 0} ]$, equally spaced in the log
scale, where $\gamma_{ X 0}$ is
calculated using the first formula in (\ref{eqinitialvalue}) with
$|Y_i - Y_j|$ replaced by $\|X_i - X_j\|$.
The rationale
for this formula can be found in
\citet{LiArtLi}.

The pair $(\gamy, \epsy)$ is selected in the same way, except that
$\epsilon_{Y 0}$ is set to 0.001. This
is because $Y$ has dimension $1$, so a weaker penalty is needed.


\section{Simulations and data analysis}\label{sec8}

In this section we present simulation comparisons among GSIR, GSAVE,
KSIR and KCCA. For the reasons explained
in the previous section, we compare GSIR with KSIR and KCCA in settings
where the sufficient predictor appears
in the conditional mean, and we compare GSAVE with GSIR, KSIR and KCCA
in settings where the sufficient
predictor appears in the conditional variance. We also apply GSIR, KSIR
and KCCA to two real data sets.

\subsection{Simulation comparisons}\label{sec8.1}

To make a comprehensive comparison of GSIR, KSIR and KCCA we consider
three regression models, namely:
\begin{eqnarray}
\cases{ \hphantom{\mathrm{II}}\mathrm{I}\mbox{:}\quad Y = \bigl(X_1^2+X_2^2
\bigr)^{1/2} \log\bigl(X_1^2+X_2^2
\bigr)^{1/2} + \varepsilon;
\vspace*{1pt}\cr
\hphantom{\mathrm{I}}\mathrm{II}\mbox{:}\quad Y = X_1/
\bigl(1+e^{X_2}\bigr)+\varepsilon,
\vspace*{1pt}\cr
\mathrm{III}\mbox{:}\quad Y = \sin\bigl(
\pi(X_1+X_2)/10\bigr)+\varepsilon, }\nonumber\\
&&\eqntext{\varepsilon\indep X,
\varepsilon\sim N(0, 0.25), p = 10;}
\end{eqnarray}
as well as three distributional scenarios for the predictor vector $X$,
namely: (A)~independent Gaussian predictors,
(B) independent non-Gaussian predictors and (C)~correlated Gaussian
predictors. In symbols:
\[
\cases{ \mathrm{A}\mbox{:}\quad X \sim N(0, I_p );
\vspace*{1pt}\cr
\mathrm{B}\mbox{:}\quad
X \sim(1/2) N(-1_p, I_p) + (1/2) N(1_p,
I_p);
\vspace*{1pt}\cr
\mathrm{C}\mbox{:}\quad X \sim N\bigl(0, 0.6 I_p + 0.4
1_p 1_p^{\mathsf{T}}\bigr). } %
\]
Note that the central $\sigma$-fields for
the three models I, II and III are generated by $X_1^2 + X_2^2$,
$X_1/(1+e^{X_2})$ and $\sin(\pi(X_1+X_2)/10)$, respectively.

We assess the quality of an estimated sufficient predictor by its
closeness to the true sufficient predictor
and its closeness to the response.
Since we are only interested in monotone functions of the predictor, we
use Spearman's correlation as the measure of closeness.
For each combination of the models and scenarios, we generate $n = 200$
observations on
$(X,Y)$ as the training data, and compute the first predicting function
using the each of three methods.
We then independently generate $m = 200$ observations on $(X,Y)$ as the
testing data, and evaluate
the predicting functions at these points.
Finally, we compute the mentioned
Spearman's correlations from the testing data.
This process is repeated $N = 200$ times. In Table~\ref{tabl1} we list means and
standard deviations of the Spearman's correlations
computed using the $N = 200$ simulated samples. From the table we see
that the performances of KCCA and GSIR are similar, and both are
slightly better
than KSIR.\vadjust{\goodbreak}

\begin{table}
\tabcolsep=0pt
\caption{Comparison of KSIR, KCCA and GSIR when sufficient predictors
appear in the conditional means}\label{tabl1}
\begin{tabular*}{\tablewidth}{@{\extracolsep{\fill}}lccccccc@{}}
\hline
\multicolumn{2}{@{}c}{\textbf{Models}}& \multicolumn{3}{c}{\textbf{Spearman cor. with
true predictor}} &\multicolumn{3}{c@{}}{\textbf{Spearman cor. with response}}
\\[-4pt]
\multicolumn{2}{@{}c}{\hrulefill}& \multicolumn{3}{c}{\hrulefill}
&\multicolumn{3}{c@{}}{\hrulefill}\\
$\bolds{X}$ & $\bolds{Y|X}$ & \textbf{KSIR} & \textbf{KCCA} & \textbf{GSIR}
& \textbf{KSIR} & \textbf{KCCA} & \textbf{GSIR}\\
\hline
A& I & 0.78 (0.05) & 0.81 (0.04)& 0.80 (0.05) &
0.63 (0.06)& 0.66 (0.05)& 0.64 (0.05) \\
& II & 0.81 (0.05) & 0.90 (0.03)& 0.91 (0.03) & 0.56 (0.06)& 0.61
(0.05)& 0.62 (0.05) \\
& III & 0.76 (0.06) & 0.89 (0.04)& 0.91 (0.03) & 0.47 (0.07)& 0.56
(0.05)& 0.56 (0.05) \\
[6pt]
B& I & 0.88 (0.02) & 0.88 (0.02)& 0.87 (0.02) &
0.82 (0.03)& 0.81 (0.03)& 0.80 (0.03) \\
& II & 0.89 (0.03) & 0.93 (0.02)& 0.93 (0.02) & 0.71 (0.04)& 0.74
(0.04)& 0.74 (0.04) \\
& III & 0.90 (0.02) & 0.97 (0.01)& 0.97 (0.01) & 0.72 (0.04)& 0.77
(0.03)& 0.77 (0.03) \\
[6pt]
C& I & 0.79 (0.04) & 0.82 (0.04)& 0.81 (0.04) &
0.64 (0.05)& 0.66 (0.05)& 0.65 (0.05) \\
& II & 0.83 (0.05) & 0.86 (0.06)& 0.88 (0.04) & 0.56 (0.06)& 0.59
(0.06)& 0.60 (0.06) \\
& III & 0.83 (0.06) & 0.96 (0.02)& 0.96 (0.02) & 0.56 (0.06)& 0.65
(0.04)& 0.65 (0.04) \\
\hline
\end{tabular*}
\end{table}

Next, we compare GSAVE, KSIR, KCCA and GSIR when the predictors only
affect the variance.
We use the following models:
\[
\cases{ \mathrm{IV}\mbox{:}\quad Y = X_1 \varepsilon;
\vspace*{1pt}\cr
\hphantom{\mathrm{I}}\mathrm{V}\mbox{:}\quad Y = (1/50) \bigl(X_1^3+X_2^3
\bigr) \varepsilon;
\vspace*{1pt}\cr
\mathrm{VI}\mbox{:}\quad Y = \bigl(X_1/
\bigl(1+e^{X_2}\bigr)\bigr) \varepsilon, } %
\]
and again the scenarios (A), (B) and (C) for the distribution of $X$.
The specifications
of $n, m, N, p$ are the same as in the previous comparison.

Because the sufficient predictors appear in the conditional
variance\break
$\var(Y|X)$ only,
it is less meaningful to measure the closeness between the estimated
sufficient predictor and the response.
So in Table~\ref{tabl2} we only report the means and standard deviations of
Spearman's correlations between the estimated and true
sufficient predictors.
We see that
GSAVE performs substantially better than the other methods. The
discrepancy can be explained by the fact
that KSIR, KCCA and GSIR depend completely on $E[\var(f(X) | Y)]$,
whereas GSAVE extracts more information from\break $\var(f(X) |Y)$.

\begin{table}[t]
\caption{Comparison of KSIR, KCCA, GSIR and GSAVE when sufficient
predictors appear in conditional variances}\label{tabl2}
\begin{tabular*}{\tablewidth}{@{\extracolsep{\fill}}lccccc@{}}
\hline
\multicolumn{2}{@{}c}{\textbf{Models}} & \multicolumn{4}{c@{}}{\textbf{Spearman's
correlation with true predictors}} \\[-4pt]
\multicolumn{2}{@{}c}{\hrulefill} & \multicolumn{4}{c@{}}{\hrulefill}
\\
$\bolds{X}$ & $\bolds{Y|X}$ & \textbf{GSAVE} & \textbf{KSIR} & \textbf{KCCA}
& \textbf{GSIR}\\
\hline
A& IV & 0.89 (0.08)& 0.10 (0.07)& 0.36 (0.22)& 0.41
(0.23) \\
& V & 0.73 (0.19)& 0.09 (0.07)& 0.17 (0.13)& 0.20 (0.14) \\
& VI & 0.84 (0.09)& 0.10 (0.08)& 0.25 (0.17)& 0.27 (0.17) \\
[6pt]
B& IV & 0.87 (0.08)& 0.10 (0.07)& 0.43 (0.25)& 0.53
(0.25) \\
& V & 0.88 (0.06)& 0.09 (0.07)& 0.11 (0.08)& 0.11 (0.08) \\
& VI & 0.76 (0.15)& 0.27 (0.11)& 0.61 (0.13)& 0.64 (0.13) \\
[6pt]
C& IV & 0.76 (0.20)& 0.11 (0.07)& 0.23 (0.16)& 0.26
(0.18) \\
& V & 0.82 (0.14)& 0.10 (0.07)& 0.11 (0.09)& 0.12 (0.09) \\
& VI & 0.73 (0.15)& 0.15 (0.10)& 0.41 (0.17)& 0.44 (0.17) \\
\hline
\end{tabular*}
\end{table}

\begin{figure}[b]

\includegraphics{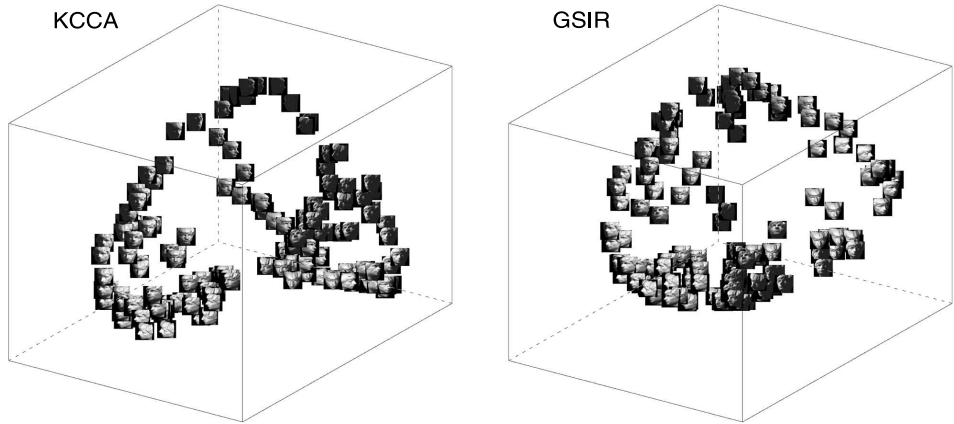}

\caption{First 3 sufficient predictors by KCCA (left panel) and GSIR
(right panel), computed from 558 training images, and evaluated on 140
testing images---faces data.}\label{figur1}\label{figureface}
\end{figure}

\subsection{Data analysis}\label{sec8.2}

We first consider the \textit{faces data}, available at
\href{http://waldron.stanford.edu/isomap/datasets.html}{http://}
\href{http://waldron.stanford.edu/isomap/datasets.html}{waldron.stanford.edu/isomap/datasets.html}. This data set
contains 698 images of the same sculpture of a face photographed at
different angles and with different lighting directions. The predictor
comprises $64 \times64$ image pixels (thus $p = 4096$), and the
response comprises horizontal rotation, vertical rotation and lighting
direction measurements (thus $q = 3$). We use this data to demonstrate
that the first three sufficient predictors estimated by KCCA and GSIR
can effectively capture the 3-variate response. We use $n = 558$ of the\vadjust{\goodbreak}
images selected at random (roughly~$80\%$) as training data, and the
remaining $m = 140$ images as testing data. For each method, we
estimate the first three predictor functions from the training data,
and evaluate them on the testing data. The left panel of Figure
\ref{figur1} is the perspective plot of the first three KCCA predictors
evaluated on the 140 testing images, and the right panel is the
counterpart for GSIR. We did not include KSIR in this comparison
because in its proposed form it cannot handle multivariate responses.
The perspective plots indicate that nearby regions in the 3-D cubes
have similar patterns of
left--right rotation, up--down rotation and lighting direction, while
distant regions have discernibly different patterns.
This reflects the ability of the three sufficient predictors to capture
the 3-variate responses.\looseness=-1

\begin{figure}

\includegraphics{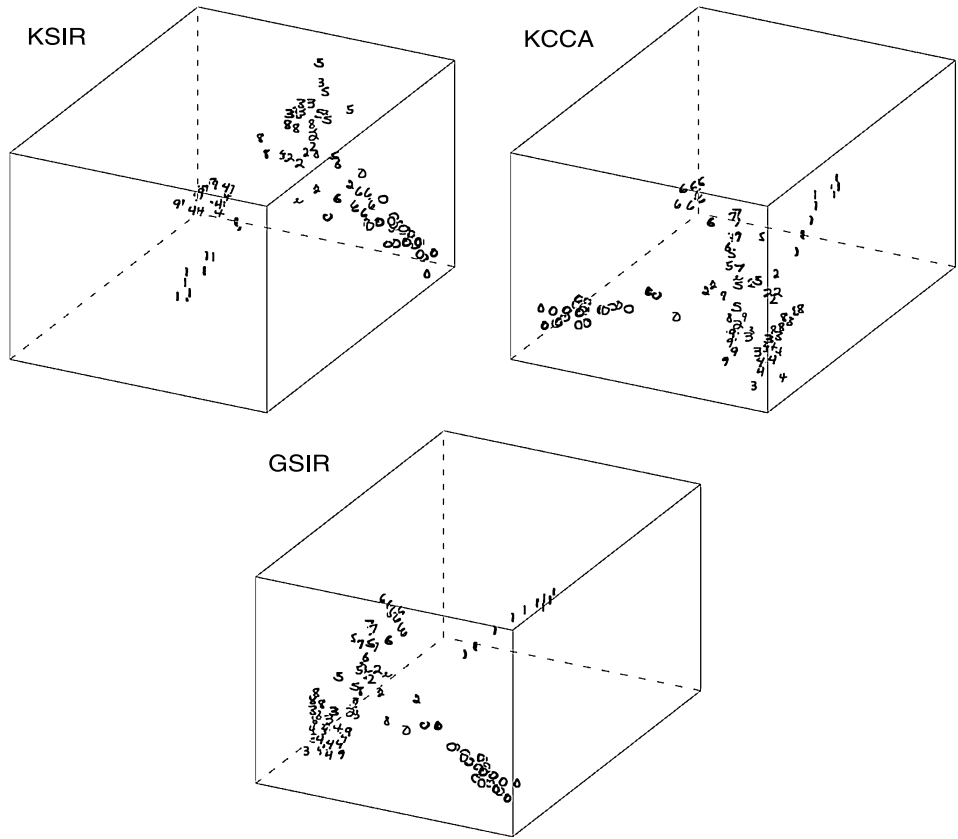}\vspace*{-3pt}

\caption{First 3 sufficient predictors by KSIR (upper-left panel), KCCA
(upper-right panel) and GSIR (lower panel), computed on 1000 training
images and evaluated on 1000 testing images---handwritten digits
data.}\label{figuredigit}\label{figur2}\vspace*{-3pt}
\end{figure}

Next, we apply KSIR, KCCA and GSIR to the \textit{handwritten digits
data}, available at
\url{http://www.cs.nyu.edu/\textasciitilde roweis/data.html}. This data set
contains 2000 images of $p = 16 \times16$ pixels showing
handwritten digits from 0 to 9---the response is thus categorical with
10 levels. We use 1000
images as training data and 1000 as testing data.
Again, for each method we estimate the first three sufficient
predictors on the
training data, and evaluate them on the testing data. Results are
presented in the three perspective
plots in Figure~\ref{figuredigit}---for visual clarity, these plots
include only 100 randomly
selected points from the 1000 in the testing data. The plots show that
all three methods provide
low-dimensional representations in which the digits are well
separated.\vspace*{-3pt}

\section{Concluding remarks}\label{sec9}

In this article we described a novel and very general theory of
sufficient dimension reduction.
This theory allowed us to combine linear\vadjust{\goodbreak} and nonlinear SDR into a
coherent system, to link them
with classical statistical sufficiency, and to subsume several existing
nonlinear SDR
methods into a unique framework.

Our developments thus revealed important and previously unexplored
properties of SDR methods.
For example, unbiasedness of various nonlinear extensions of SIR
proposed in recent literature was
proved under the stringent linear conditional mean assumption. We were
able to show that these
methods are all unbiased under virtually no assumption, and that GSIR
is exhaustive under the completeness
assumption. We were also able to show that nonlinear extensions of SIR
are in general \textit{not} exhaustive
when completeness is not satisfied, and that in these cases GSAVE can
recover a larger portion of the central
class. These
insights could not have been obtained without paralleling linear and
nonlinear SDR as
allowed by our new theory.\looseness=-1

In addition to achieving theoretical synthesis and important insights
on SDR methods, we introduced a
new \textit{heteroscedastic} conditional variance operator---which is
more general than the (homoscedastic) conditional
variance operator in Fukumizu, Bach and Jordan (\citeyear{FukBacJor03}, \citeyear{FukBacJor09}). This
operator was crucial to generalizing SAVE to the
nonlinear GSAVE, and thus to exploit dependence information in the
conditional variance to improve upon the
performance of the nonlinear extensions of SIR.
We have no doubt that the heteroschedastic conditional variance
operator can be used to generate nonlinear
extensions of other second-order SRD methods such as
contour regression [\citet{LiZhaChi05}], directional
regression [\citet{LiWan07}], SIR-II [\citet{Li91N2}]
and other F2M methods [\citet{CooFor09}]. These extensions will
be the topic of future work.

More generally, it is our hope that the clarity and simplicity that
classical notions lend to the formulation of dimension
reduction, as well as the transparent parallels we were able to draw
between linear and nonlinear SDR, will provide
fertile grounds for much research to come.

As we put forward a general theory that encompasses both linear and nonlinear
SDR, it is also important to point out that linear SDR
has its special values that cannot be replaced by nonlinear SDR via
kernel mapping, one of which is
its preservation of the original coordinates and as a result its strong
interpretability.
For example, when mapped to higher dimension spaces, kernel methods
can sometimes interpret difference in variances in the original
coordinates as location separation in the transformed
coordinates, which can
be undesirable depending on the goal and emphasis of particular
applications. For further discussion
and an example of this point, see \citet{LiArtLi}.


\section*{Acknowledgments}

We would like to thank two referees and an Associate Editor for their
insightful comments and useful suggestions, which led to significant
improvements to this paper. In particular, the consideration of
nonlinear sufficient dimension reduction for functional data is
inspired by the comments of two referees.\vadjust{\goodbreak}

\begin{supplement}
\stitle{Supplement to ``A general theory for nonlinear sufficient
dimension reduction: Formulation and estimation''}
\slink[doi]{10.1214/12-AOS1071SUPP} 
\sdatatype{.pdf}
\sfilename{aos1071\_supp.pdf}
\sdescription{This is supplementary appendix that contains
some techincal proofs of the results in the paper.}
\end{supplement}


\printaddresses

\end{document}